\documentstyle[11pt]{article}
 \normalbaselines

\font\ggordo=cmbx10 scaled \magstep2

\font\ninerm=cmr9
\font\nineit=cmti9
\font\ninebf=cmbx9

 \newtheorem{definition}{Definition}
 \newtheorem{theorem}{Theorem}
 
 \newtheorem{lemma}{Lemma}
 \newtheorem{proposition}{Proposition}
 \newtheorem{remark}{Remark}

 \def \R{{\rm I\kern -2.2pt R\hskip 1pt}}
 \def \Z{{\rm Z\!\!Z}}
 \def \Q{\hskip2pt {\rm 0\kern -7pt Q}}
 \def \P{{\rm
 I\kern -2.2pt P\hskip 1pt}}

 \def \N{{\rm I\kern -2.1pt N\hskip 1pt}}
 \def \C{{\rm C\kern -4.6pt
      \vrule height 6.8pt width 0.3pt depth -0.5pt}\hskip4pt}
 \def \K {{\rm
 I\kern -2.1pt K\hskip 1pt}}

 \def\cqfd{\vbox{\hrule height 5pt  width 5pt }\bigskip}
 \def\noi{\noindent}

\begin{document}

\title{\ggordo STRAIGHT--LINE PROGRAMS IN GEOMETRIC ELIMINATION THEORY}

 \author{ M. Giusti $^{1,3}$ \
 J. Heintz $^{2,3}$
 J.E. Morais$^{2,3}$ J. Morgenstern $^{}$ L.M. Pardo $^{2,3}$ }%

\date{}

 \maketitle

\addtocounter{footnote}{1}\footnotetext{\nineit : GAGE, Centre de
Math\'ematiques. \'Ecole Polytechnique. F-91228, Palaiseau Cedex.
\\ FRANCE}

\addtocounter{footnote}{1}\footnotetext{\nineit : Dept. de Matem\'aticas,
Estad\'{\i}stica y Computaci\'on. F. de Ciencias. U. Cantabria. \\ E-39071 
SANTANDER, Spain}

\addtocounter{footnote}{1}\footnotetext{\nineit : Research was partially
supported by the following European, French and Spanish grants~: \\ PoSSo
ESPRIT/BRA 6846,  GDR CNRS 1026 MEDICIS, DGCYT PB92--0498--C02--01 and
\\ PB93--0472--C02--02}

\begin{center} {\bf Dedicated to Volker Strassen for his work on complexity}
\end{center}

\begin{abstract}
 We present a new method for solving symbolically zero--dimensional polynomial
equation systems in the affine and toric case. The main feature of  our method
is the use of problem adapted data structures~: arithmetic networks  and
straight--line programs. For sequential time complexity  measured by network
size we obtain  the following result~: it is possible to solve any affine or
toric zero--dimensional equation system in non--uniform sequential time which is
polynomial in the length of the input description and the ``geometric degree" of
the equation system. Here, the  input is thought to be given by a straight--line
program (or alternatively in  sparse representation), and the length of the
input is measured by number of  variables, degree of equations and size of the
program (or sparsity of the  equations). The geometric degree of the input
system has to be adequately defined. It is  always bounded by the
algebraic--combinatoric ``B\'ezout number" of the system  which is given by the
Hilbert function of a suitable homogeneous ideal. However, in many
 important cases, the value of the geometric degree of the system is much
smaller than its B\'ezout number since this geometric degree does not take into
account multiplicities or degrees of extraneous components (which may appear at
infinity in the affine case or may be contained in  some coordinate hyperplane
in the toric case).

Our method contains a new application of a classic tool to symbolic
computation~: we use Newton iteration in order to simplify straight--line
programs occurring in elimination procedures. Our new technique allows for
practical implementations a meaningful characterization of the intrinsic {\it
algebraic complexity} of typic elimination problems and reduces the still
unanswered question of their intrinsic {\it bit complexity} to algorithmic
arithmetics. However our algorithms are not rational anymore as are the usual
ones in elimination theory. They require some restricted computing with
algebraic numbers. This is due to its numeric ingredients (Newton iteration).
Nevertheless, at least in the case of polynomial equation systems depending on
parameters, the practical advantage of our method with respect to more
traditional ones in symbolic and numeric computation is clearly visible.  Our
approach produces immediately a series of division theorems (effective
Nullstellens\"atze) with new and more differentiated degree and complexity
bounds (we shall state two of them).

It should be well understood that our method does not improve the well known 
worst--case complexity bounds for zero--dimensional equation solving in
 symbolic and numeric computing.

Part of the results of this paper were announced in \cite{gihemopar}.
\end{abstract}

\smallskip

\noi {\bf Keywords.} {\small Polynomial system solving, elimination, division
theorem, Nullstellensatz, geometric  degree, arithmetic network, straight--line
program, complexity.}

\section{Introduction and Statement of Results}
\label{introduccion} Let $k$ be an infinite, perfect field which we think 
``effective" with respect to arithmetic operations as addition/subtraction,
multiplication/division and extraction of $p$--th roots in case that $k$ has
positive  characteristic $p$. Let ${\bar k}$ be a fixed algebraic closure of
$k$.

An important problem in elimination theory is the computation of the isolated
points of an affine algebraic variety. This means we consider the algebraic
variety
$$V\!:=V(f_1,\ldots,f_s)\!:=\{ x\in {\bf A}^n({\bar k}) : \;
f_1(x)=0,\ldots,f_s(x)=0
\},$$  where $f_1,\ldots,f_s\in k[X_1,\ldots,X_n]$ with $ s \ge n$ are
polynomials of degree at most $d$ in the variables $X_1,\ldots,X_n$ and we look
at the  following problem~: given a nonzero linear form
$H \in  k[X_1,\ldots,X_n]$ and a new variable $T$, compute a nonzero polynomial
$p\in k[T]$  such that $p(H)$ vanishes on all isolated points of $V$.  This
``weak" form of symbolic solving of the equation system given by the polynomials
$f_1,\ldots,f_s$ is in all thinkable aspects (effectivity, practical efficiency
and theoretical complexity) equivalent to  the following more explicit form~:
find {\it univariate} polynomials
$q,v_1,\ldots,v_n \in k[T]$ such that  the set of isolated points of $V$ can be
written as $\{(v_1(t),\ldots,v_n(t))\; ;\; t \in {\bf A}^1({\bar k}) \;,\;
q(t)=0\}$ (see \cite{giu-he92},
\cite{krick-pardo1} and Lemma \ref{primitive} below).

As many authors do (see for instance, \cite{ber-yg1}, \cite{ber-yg2}, 
\cite{brown1},
\cite{can-gal-he}, \cite{cagahe}, \cite{can}, \cite{chigri}, \cite{dube}, 
\cite{figismi},
\cite{giu-he92}, \cite{gihesa}, \cite{kollar}, \cite{krick-pardo1},
\cite{krick-pardo:CRAS},
\cite{lala}, \cite{phil}, \cite{saso}), we replace the original input system by 
$n$ generic $k$--linear  combinations of the equations $f_1,\ldots,f_s$. This
preparation of the input conserves all irreducible components (and, in
particular, the isolated points) of $V$ and adds possibly some new, extraneous
ones (these components have to be eliminated afterwards). Thus, let us suppose
from now on that $s=n$ holds. This implies that the input system
$f_1,\ldots,f_n$ forms locally a complete intersection with respect to the
isolated points of $V$. In view of the surprisingly low (linear) regularity
bound for homogeneous global complete intersection ideals (see \cite{bria},
\cite{laz}, \cite{macau}) one is tempted to replace the original equations 
$f_1,\ldots,f_n\in k[X_1,\ldots,X_n]$ by new homogeneous ones,
$G_1,\ldots,G_n
\in k[X_0,\ldots,X_n]$, and the affine variety
$V=V(f_1,\ldots,f_n)$ by the corresponding projective one, namely  
$W\!:=V(G_1,\ldots,G_n)=\{x \in \P^n({\bar k}) \; ; \;
G_1(x)=0,\ldots,G_n(x)=0\}$. A simple minded way to proceed consists in taking
as $G_1,\ldots,G_n$ just the homogenizations of the original system
$f_1,\ldots,f_n$. However this idea does not produce the desired effect even in
 case that the polynomials
$f_1,\ldots,f_n$ form a regular sequence in $ k[X_1,\ldots,X_n]$. The reason  is
that the homogenizations of $f_1,\ldots,f_n$ need not to generate a complete
intersection ideal.

In order to remedy this defect, many authors use homotopic deformations and
path--following methods for the construction of a regular sequence of auxiliar
homogeneous polynomials $G_1,\ldots,G_n$, from whose zeroes the solutions of the
original system $f_1,\ldots,f_n$ are then extracted (see
\cite{royetal}, \cite{can}, \cite{chis}, \cite{chigri}, \cite{dicketal},
\cite{gar-zan}, \cite{giu-he92}, \cite{lalo-lupe}, 
\cite{krick-pardo1}, 
\cite{krick-pardo:CRAS}, \cite{lala}, \cite{moeller}, 
\cite{shub-sma1}). The polynomials $G_1,\ldots,G_n$ are homogeneous in the
variables
$X_0,X_1,\ldots,X_n$. Moreover they depend on a deformation parameter
$\varepsilon$ and form a regular sequence in 
$k(\varepsilon)[X_0,X_1,\ldots,X_n]$. They define a zero dimensional projective
subvariety of $\P^n(\overline {k(\varepsilon)})$ without points at infinity
(here $\overline {k(\varepsilon)}$ denotes any algebraic closure of the function
field $k(\varepsilon)$). Specializing the polynomials
$G_1,\ldots,G_n$ in $\varepsilon=0$ we obtain a projective subvariety of
$\P^n({\bar k})$ which contains the isolated points of $V$ as irreducible
components. Since the underlying deformation is flat in the isolated points of
$V$, well known techniques based on the implicit or explicit use of Macaulay's
{\it u}--resultant allow first to find the solutions of
$V(G_1,\ldots,G_n)$ and finally the isolated points of $V$.

This deformation method introduces a somehow artificial dependency of the
complexity of the algorithms on the regularity of the Hilbert function of the
homogeneous ideal $(G_1,\ldots,G_n)$ of $k(\varepsilon)[X_0,X_1,\ldots,X_n]$ and
on the degree of $(G_1,\ldots,G_n)$ which is defined by means of the Hilbert
polynomial and is called its B\'ezout number. The regularity of
$(G_1,\ldots,G_n)$ is bounded by $nd-n$. This implies that deformation based
algorithms have to triangulate rectangular matrices of size 
$\displaystyle{ { {nd}
\choose {n}}} \times n\displaystyle{ { {nd}
\choose {n}}}$ in order to extract from this data some ``geometrically
meaningful" square matrix of size the B\'ezout number. The characteristic
polynomial of this square matrix is the basic eliminating form that all
algorithms work with. However, the B\'ezout number of the ideal
$(G_1,\ldots,G_n)$ (i.e. the degree of the eliminating form) includes
multiplicities and degrees of extraneous components which the previous
deformation process adds to the original variety $V$. The exact value of the
B\'ezout number is $\prod\limits_{i=1}^n deg G_i$ which is of order $d^n$.

The method we present in this paper will be independent of  ``algebraic"
quantities such as regularity of the Hilbert function or B\'ezout number of an
appropriate homogeneous ideal. This allows in geometrically well suited cases
(typically when the geometric degree of the variety
$V$ is low) to reduce the size of the matrices in the algorithms (and hence the
complexity of the procedures). However, in worst case (e.g. when the equations
$f_1,\ldots,f_n$ are generic) our complexity bounds are roughly the same as
those obtained by deformation based algorithms. This is not surprising in view
of the known lower bounds for elimination problems (
\cite{he-morg},
\cite{shub-sma2}).

An important aspect of elimination procedures consists in the encoding of the
polynomials appearing as inputs, outputs or intermediate results of the
algorithm. Encoding polynomials by their coefficients (dense representation)
faces us with an input size of order
$\displaystyle{{{d+n}
\choose n}\leq c \cdot d^n}$ (if $n,d\ge 3$ we can take $c=1$, otherwise we can
take as $c$ the base of the natural logarithm). If the equations
$f_1,\ldots,f_n$ are generic, the geometric degree of the variety
$V$ and the B\'ezout number  of the system (given by the homogenizations of
$f_1,\ldots,f_n$) coincide and are of order
$d^n$. Thus, a complexity bound of type $d^{O(n)}$ is both polynomial in the
dense input size and the degree of the variety $V$, what is the best we can hope
in worst case (see \cite{he-morg} and \cite{shub-sma2}). 

However, when the geometric degree of $V$ is low, one wishes to use more
economic encodings of inputs, outputs and intermediate results. Such an encoding
is given, for instance, by the data structure straight--line program.

For the moment let us fix the following notations and assumptions~:

There is given a family of $n$ input polynomials $f_1,\ldots,f_n \in
k[X_1,\ldots,X_n]$ which are thought to be encoded by a straight--line program
(arithmetic circuit) $\beta$ without essential divisions in
$k[X_1,\ldots,X_n]$ (this means we allow $\beta$ only to contain divisions by
non--zero elements of $k$). We denote the nonscalar size (``nonscalar length" in
more traditional terminology) and the nonscalar depth of
$\beta$ by
$L$ and
$\ell$ respectively (see \cite{gathen:par}, \cite{joos:aaecc5},
\cite{krick-pardo1}, \cite{Strassen:hand} and Section \ref{modelo} for the
notions of straight--line program, arithmetic network and the complexity
measures we shall use subsequently).

Let us first consider the {\it affine case}. Here, we suppose that
$f_1,\ldots,f_n$ form a regular sequence in $k[X_1,\ldots,X_n]$. Fix $1
\le i
\le n$. The affine variety defined by the ideal
$(f_1,\ldots,f_i)$, namely $V(f_1,\ldots,f_i)$, is a Zariski closed subset of
${\bf A}^n\!:= {\bf A}^n({\bar k})$ of {\it pure dimension} $n-i$ (i.e. all 
irreducible components of $V(f_1,\ldots,f_i)$ have the same dimension
$n-i$; we say also that
$V(f_1,\ldots,f_i)$ is an {\it equidimensional variety of dimension}
$n-i$). The {\it (geometric) affine degree} of $V(f_1,\ldots,f_i)$ is defined as
usual as the cardinality of the finite set of points we obtain  cutting
$V(f_1,\ldots,f_i)$ by $n-i$ generic affine hyperplanes of ${\bf A}^n$ (more
generally, we define the affine degree of a closed Zariski subset of
${\bf A}^n$ as the sum of the degrees of its irreducible components. See e.g.
\cite{joos:tesis} or
\cite{fulton} for this notion of degree and its motivation). We denote the
affine degree of $V(f_1,\ldots,f_i)$ by $deg V(f_1,\ldots,f_i)$. 

Let us now define $\delta \!:= max \{deg V(f_1,\ldots,f_i)\; ; \; 1 \le i
\le n
\}$ as the {\it (geometric) affine degree of the equation system}
$f_1,\ldots,f_n$. We write $V\!:= V(f_1,\ldots,f_n)$ and observe that $V$
contains only finitely many, namely $degV \le \delta$ points.

Let us now consider the {\it toric case}. Fix $1 \le i \le n$. We define the
{\it toric irreducible components} of the affine variety
$V(f_1,\ldots,f_i)$ as those which are not contained in one of the hyperplanes
$V(X_j)$, $1 \le j \le n$, of
${\bf A}^n$. We suppose that $f_1,\ldots,f_i$ form a {\it toric complete
intersection}. This means that toric irreducible components of
$V(f_1,\ldots,f_i)$ exist having all dimension  $n-i$. The {\it (geometric)
toric degree} of
$V(f_1,\ldots,f_i)$ is denoted by $deg^{*}V(f_1,\ldots,f_i)$ and is defined as
the sum of the affine degrees of all toric components of
$V(f_1,\ldots,f_i)$. We call $\delta^{*}\!:= max \{deg^{*} V(f_1,\ldots,f_i)\; ;
\; 1 \le i \le n\}$  the {\it (geometric) toric degree of the equation system}
$f_1,\ldots,f_n$.

\medskip
 
{\it Let be given a non--zero linear form $H$ of $k[X_1,\ldots,X_n]$ represented
by its coefficients and let $T$ be a new variable. We consider as the fixed
input of all our principal algorithms the division free straight--line program
$\beta$ which computes $f_1,\ldots,f_n$ and, where it makes sense, the
coefficient $n$--tuple of $H$. We recall that
$d$ is an upper bound for the degrees of the polynomials of the input system
$f_1,\ldots,f_n$, that $\delta$ and $\delta^{*}$ represent the geometric affine
degree  and the geometric toric degree of the system
$f_1,\ldots,f_n$ respectively. Furthermore we recall that $L$ and $\ell$ are
the nonscalar size and depth of $\beta$ respectively}.  

With these notations and assumptions, we may state our main results as
follows~:

\begin{theorem}[the affine case]
\label{afin}  Assume that $f_1,\ldots,f_n$ form a regular sequence in

\noindent $k[X_1,\ldots,X_n]$. Suppose furthermore that for $1Ê\le i \le n$ the
ideal $(f_1,\ldots,f_i)$ generated by $f_1,\ldots,f_i$ in $k[X_1,\ldots,X_n]$ is
radical. Then there exists an  arithmetic network with parameters in ${\bar k}$
which, from the input given by
$\beta$ and
$H$, computes the coefficients of a monic polynomial $p\in k[T]$ such  that
$p(H)$ vanishes on
$V=V(f_1,\ldots,f_n)$.The nonscalar size and depth of the network are 
$(nd\delta L)^{O(1)}$ and
$O(n(log_2(n d) +\ell)log_2 \delta)$ respectively.
\end{theorem}

\begin{theorem}[the toric case]
\label{tor} Assume that for $ 1 \le i \le n$ the polynomials
$f_1,\ldots,f_i$ form a toric complete intersection and generate an ideal in
$k[X_1,\ldots,X_n]$ whose localization
$(f_1,\ldots,f_i)_{\prod\limits_{j=1}^n X_j}$ is radical. Then there exists an 
arithmetic network with parameters in ${\bar k}$ which, from the input given by
$\beta$ and
$H$, computes the coefficients of a monic polynomial $p^{*}\in k[T]$ such  that
$p^{*}(H)$ vanishes on the toric variety $V^{*}\!:= V \setminus
V(\prod\limits_{j=1}^n X_j)$. The nonscalar size and depth of the network are 
$(nd\delta^{*} L)^{O(1)}$ and
$O(n(log_2(n d) +\ell)log_2 \delta^{*})$  respectively.
\end{theorem}

Let us remark that the algorithms underlying Theorems \ref{afin} and
\ref{tor} make substantial use of linear algebra subroutines dealing with
matrices of size at most $2d \delta$ (affine case) and $2d
\delta^{*}$ (toric case).

We may sharpen the assertions of Theorem \ref{afin} and \ref{tor} to the
following statement~:

\begin{theorem}
\label{parametrizacion} Let $f_1,\ldots,f_n$ be polynomials of
$k[X_1,\ldots,X_n]$ which satisfy the assumptions of Theorems \ref{afin} or
\ref{tor}. Then there exists an arithmetic network with parameters in
${\bar k}$ which from the input given by the circuit $\beta$ computes in the
affine case the coefficients of polynomials
$q,v_1,\ldots,v_n \in k[T]$ and in the toric case  the coefficients of
polynomials $q^{*},v^{*}_1,\ldots,v^{*}_n \in k[T]$ such that the following
conditions hold~:

\begin{enumerate}
\item  $q$ and $q^{*}$ are monic and have degrees $degq=deg V$
      and $degq^{*}=degV^{*}=deg^{*}V$ respectively. Moreover
      $v_1,\ldots,v_n$
      and $v^{*}_1,\ldots,v^{*}_n$ satisfy the degree bounds $max \{deg v_i
      \; ;\; 1 \le i\le n\}< degq$ and $max \{deg v^{*}_i \; ;\; 1 \le i\le
      n\}< degq^{*}$.

\item the zerodimensional varieties $V$ and $V^{*}$ can be  
      parameterized as 
       
      $$V\!=\{(v_1(t),\ldots,v_n(t))\; ;\; t \in {\bf
      A}^1\;,\; q(t)=0\} \; and$$ $$V^{*}\!=\{(v^{*}_1(t),\ldots,v^{*}_n(t))\;
;\;
      t \in {\bf A}^1 \;,\; q^{*}(t)=0\}$$.
\end{enumerate} The nonscalar size and depth of the network are $(nd\delta
L)^{O(1)}$ and
$O(n(log_2(n d) +\ell)log_2 \delta)$ in the affine case and 
$(nd\delta^{*} L)^{O(1)}$ and
$O(n(log_2(n d) +\ell)log_2 \delta^{*})$ in the toric case.
\end{theorem}

Of course we can reformulate Theorems \ref{afin}, \ref{tor} and
\ref{parametrizacion} for sparse inputs. In this case the parameters $L$ and
$\ell$ measuring nonscalar size and depth of the input circuit $\beta$ can be
replaced by the sparsity $N$ of the polynomials $f_1,\ldots,f_n$ and by
$log_2d$ respectively (the sparsity is the number of monomials with non--zero
coefficients appearing in $f_1,\ldots,f_n$). The nonscalar size and depth of the
resulting network are then $(nd \delta N)^{O(1)}$ and
$O(n log_2(nd) log_2
\delta)$ in the affine case and $(nd \delta^{*} N)^{O(1)}$ and $O(n log_2(nd)
log_2 \delta^{*})$ in the toric case

The most important context of applications of  Theorems \ref{afin},
\ref{tor} and \ref{parametrizacion}  is the following situation of
``parametric/numeric equation solving"~:

{\it let $\Omega$ be a prime field with algebraic closure ${\bar
\Omega}$ and let $\theta_1,\ldots,\theta_m$ be  indeterminates over
$\Omega$ (called ``parameters"). Let
$k\!:=\Omega(\theta_1,\ldots,\theta_m)$ (there will be no harm to our arguments
in characteristic $p$ since
$\Omega$ is perfect). Suppose now that $\beta$ is a division free straight--line
program in
$\Omega[\theta_1,\ldots,\theta_m,X_1,\ldots,X_n]$ of nonscalar length and depth
$L$ and $\ell$ respectively. Thus our polynomials
$f_1,\ldots,f_n$ belong to the polynomial ring
$\Omega[\theta_1,\ldots,\theta_m,X_1,\ldots,X_n]$. Assume that $H$ is a nonzero
linear form of $\Omega[X_1,\ldots,X_n]$. The varieties $V$ and
$V^{*}$ are interpreted as subvarieties of the affine space ${\bf A}^n({\bar
k})={\bf A}^n (\overline{\Omega(\theta_1,\ldots,\theta_m)})$. The remaining
assumptions and notations are the same as in Theorems
\ref{afin}, \ref{tor} and
\ref{parametrizacion} for $k:=\Omega(\theta_1,\ldots,\theta_m)$}.

In this situation the statement of Theorem \ref{afin}, \ref{tor} and
\ref{parametrizacion} can be sharpened as follows~: 

\begin{theorem} 
\label{aplicacion}
There exists an arithmetic network with parameters from the field 
${\bar \Omega}$, which has nonscalar length $(nd \delta L)^{O(1)}$ and depth
$O(n(log_2(nd)+\ell)log_2 \delta)$ in the affine and which has nonscalar length
$(nd \delta^{*} L)^{O(1)}$ and depth $O(n(log_2(nd)+\ell)log_2 \delta^{*})$ in
the toric case. This arithmetic network produces as output a division free
straight--line program in ${\bar \Omega}[\theta_1,\ldots,\theta_m]$, say 
$\gamma$, in the affine and $\gamma^{*}$ in the toric case. This
straight--line program has the following properties~:

\begin{itemize}
\item according to the case (affine or toric) the
      circuits $\gamma$ and $\gamma^{*}$
      represent the coefficients with respect to the variable $T$ of
      polynomials $p$ and $p^{*}$ or
      $q,w_1,\ldots,w_n$ and
      $q^{*},w^{*}_1,\ldots,w^{*}_n$ which belong to
      $\Omega[\theta_1,\ldots,\theta_m,T]$. Moreover in the 
      situation where it makes sense $\gamma$ and $\gamma^{*}$ compute also the
      (nonzero) discriminants $\rho \in \Omega[\theta_1,\ldots,\theta_m]$
      and $\rho^{*} \in \Omega[\theta_1,\ldots,\theta_m]$ of $q$ and
      $q^{*}$. 
\item the straight--line programs $\gamma$ and $\gamma^{*}$ have
      asymptotically the same
      nonscalar size and depth as the arithmetic network generating them. 
\item the (monic in $T$) polynomials
      $p,p^{*},q,q^{*} \in \Omega[\theta_1,\ldots,\theta_m,T]$ and the
      rational functions  $v_1\!:={w_1 \over \rho},\ldots,{w_n \over \rho}
      \in \Omega(\theta_1,\ldots,\theta_m)[T]$ and $v^{*}_1\!:={w^{*}_1
\over
      \rho},\ldots,{w^{*}_n \over \rho}
      \in \Omega(\theta_1,\ldots,\theta_m)[T]$ have the properties
      stated
      for them in Theorems \ref{afin}, \ref{tor} and
      \ref{parametrizacion} respectively. For instance in the situation of
      Theorem
      \ref{afin} the polynomial $p(\theta_1,\ldots,\theta_m,H)$ vanishes on
      $V=V(f_1,\ldots,f_n)= \{x \in
      \overline{\Omega(\theta_1,\ldots,\theta_m)}^n
       \; ;\; f_1(x)=0,\ldots,f_n(x)=0 \}$ and in the situation of Theorem
      \ref{parametrizacion}(i) the algebraic variety $V$ has the form
      $V= \{v_1(t),\ldots,v_n(t) \; ; \; t \in
     \overline{\Omega(\theta_1,\ldots,\theta_m)} \; , \;
     q(\theta_1,\ldots,\theta_m,t)=0 \}$.
\end{itemize} 
\end{theorem}

We shall omit the proof of Theorem \ref{aplicacion} because it would be almost
textually the same as that of Proposition \ref{inter.comp.} below.

A natural question to ask is the following one~: what ``real life" complexities
are hidden  behind the notions of arithmetic network and straight--line program
with parameters in ${\bar k}$? What is the link of these data structures to
``ordinary" rational arithmetic networks and straight--line programs with
parameters in the coefficient field
$k$? For answering this question let us restrict ourselves to the affine case
(the toric case can be discussed analogously).

A simple minded translation of the arithmetic network with parameters in
${\bar k}$ underlying Theorem \ref{afin} to the rational context over $k$ would
produce an ordinary arithmetic network with parameters in $k$ of nonscalar size
and depth $(nd\delta^{n} L)^{O(1)}$ and
$O(n(log_2(d\delta) +\ell))$ respectively. Whereas the nonscalar depth is fine,
the nonscalar size of the network grows up to order of the sequential time
complexity of usual Groebner basis computations for the elimination task we are
considering, i.e. to
$(nd^{n^2})^{O(1)}$ in worst case (here we assume that $\delta=d^n$ and 
$L=\displaystyle{{{d+n} \choose n}}= O(d^n)$ holds; see \cite{can-gal-he} and
\cite{dicketal}).

The problem arising in this way can be localized in the iterative character of
the algorithm ($n$ iterations) which in the ordinary rational arithmetic network
version produces stupidly growing straight--line programs for the representation
of intermediate results. This is due to the repeated use of interpolation
subroutines which nevertheless can be summarized by suitable FOR commands (see
\cite{gihemopar}). However, applying ``data compression" by means of the already
mentioned homotopic deformation method one can lower this ordinary network size
to
$(nd^n)^{O(1)}$ which is an already known complexity bound for the elimination
task we are considering (\cite{can}, \cite{chis},
\cite{giu-he92},
\cite{krick-pardo1}, \cite{lala}. See also \cite{ber-yg1}, 
\cite{chigri}, \cite{he-morg}, \cite{laz}). This means the advantage of the
rational version of our method is reduced to exegetic ``practical complexity" as
far as running time is concerned (nevertheless we may have saved something with
respect to storage space).

However a more radical compression of the enormously expanding straight--line
programs during the execution of the algorithms is possible thanks to a symbolic
adaptation of Newton's method (or Hensel iteration as someones like to call it
in case of positive characteristic) in combination with ``Vermeidung von
Divisionen" (\cite{stras:verm}). This compression is done by the algorithm
underlying the fundamental technical Lemmas \ref{sinfor} and
\ref{sinfortorico} below.

The key point of our method consists in a consequent use of problem adapted data
structures~: arithmetic networks and straight--line programs with parameters in
${\bar k}$ (or $k$). A similar ground idea may be found in the impressive
theoretical and practical work of D. Duval and her school on dynamic evaluation
(\cite{duv}). Closest to our complexity results comes the numeric method of M.
Shub and S. Smale for finding ``approximate" solutions of zero--dimensional
homogeneous ``average" equations systems (\cite{shub-sma1}). Their method too is
independent on the B\'ezout number and the regularity of the input ideal. The
sequential time complexity of their algorithm (measured as ours by counting just
arithmetic operations) is polynomial in $n
\displaystyle{{{d+n}
\choose n}}$ (the size of the densely given input system) and the number of
approximate solutions wanted. However caution is necessary in the interpretation
of their result as ``equation solving" in the sense of computer algebra. First
they need  sufficient genericity of the system (just the contrary of the point
of view of computer algebra which focuses on ``special systems"). Secondly
Newton iteration (for which they seek approximate solutions) is inefficient in
terms of bit--complexity if the rational numbers appearing in it are given in
binary coding. This fact is due to Liouville's estimate (see
\cite{pardo}). However this second drawback can be avoided coding rational
numbers by straight--line programs. Identity of rational numbers given by such
encodings can be detected in random polynomial time by a BPP--test 
(\cite{schonhage}, \cite{iba-mor}) or a fixed sample test (as in
\cite{heschnorr} or
\cite{cucketal}) which is due to \cite{ha-mon}. Summarizing this we can say that
in case $\Omega\!:= \Q$ and $k\!:=\Q(\theta_1,\ldots,\theta_m)$ with
$\theta_1,\ldots,\theta_m$ algebraically independent, our algorithms have a
reasonable translation to the (probabilistic) bit--model and that this
translation conserves the overall complexity character of our algorithms (this
is work in progress and will be the subject of a further publication; see
Section \ref{conclusiones}). Of course it would also be advantageous to dispose
of algorithms of the same complexity type which are realizable exclusively by
{\it rational} arithmetic networks and straight--line programs (with parameters
in $k$). This is a question which we shall also consider in further work.

Our paper is organized as follows~: 

\noindent in Section \ref{modelo}  we introduce our computational model of
arithmetic networks and straight--line programs with parameters in ${\bar k}$
(or $k$). In Section
\ref{compresion} we prove Lemma \ref{sinfor} and \ref{sinfortorico} which
represent a key tool for our algorithms. Thanks to these lemmas we are able to
compress efficiently certain straight--line programs which appear as
intermediate computations in our procedures. The algorithm underlying Lemma
\ref{sinfor} and \ref{sinfortorico} requires the use of algebraic numbers
(observe that this situation produces the side  effect that our subsequent
arithmetic networks and straight--line programs depend on algebraic parameters).
In Section
\ref{lemas} we review some known technical lemmas which are necessary for the
proofs in Section
\ref{alg}. The essential steps of the proof of Theorem \ref{afin},
\ref{tor} and
\ref{parametrizacion} are contained in Section \ref{alg}. This proof is based 
on a recursive construction of a Noether normalization and a primitive element
for the varieties
$V(f_1,\ldots,f_i)$, $1
\le i \le n$, followed by a procedure of cleaning extraneous irreducible
components. This recursive construction is explained in Section \ref{alg} and a
recursion law of the complexity of the underlying algorithm is given. At this
point the application of Lemma \ref{sinfor} and \ref{sinfortorico} of Section
\ref{compresion} becomes  crucial. A simple minded iteration of the recursive
construction we introduce in Section \ref{alg} would lead to a size explosion of
straight--line programs. Lemma \ref{sinfor} and \ref{sinfortorico} allow us to
compress after each iteration step the straight--line programs occuring in our
algorithm. The core of the proof of Theorem
\ref{afin},
\ref{tor} and \ref{parametrizacion} is contained in Proposition
\ref{inter.comp.} and \ref{hipersurf} which generalizes the statement of Theorem
\ref{parametrizacion} to the case of arbitrary reduced complete intersection
ideals. This proposition is used in Section \ref{division} in order to formulate
two division theorems (effective Nullstellens\"atze) with new
complexity and degree bounds which are polynomial in our parameters $n, d,
\delta$ and
$L$ (Theorem \ref{division1} and \ref{division2}). These new division theorems
do not improve the degree (and height) bounds of the classical ones
(\cite{ber-yg1}, \cite{brown1}, \cite{can-gal-he}, \cite{cagahe}, \cite{dube},
\cite{kollar}, \cite{krick-pardo1}, \cite{krick-pardo:CRAS},  \cite{phil},
\cite{saso}) but they have three advantages~: the first one is that they are
better suited for computational issues, the second one is that they are more
general and the third one is that they explain better the classical results.

\section{The Computational Model}  
\label{modelo} The idea of using straight--line programs as succinct encodings
of special polynomials (as they appear in elimination problems) goes back to the
late seventies and was discovered by different people independently (one of them
being our coauthor J. Morgenstern who prefered oral ``hadise" to hardcover
publications for the dissemination of his thought). This idea appears implicitly
or explicitly in the following representative (although not exhaustive) list of
papers, mostly dedicated to probabilistic testing of polynomial identities~:
\cite{heschnorr}, \cite{iba-mor}, \cite{schw}, \cite{demillolipton},
\cite{zip}. First applications of this point of  view to computer algebra dealt
only with elimination of just one variable (absolute primality testing, greatest
common divisor computation and factorization of multivariate polynomials; see
\cite{hesie},
\cite{kalt1},
\cite{kalt2}) and were later extended to multivariate elimination problems by
means of ``duality techniques" of different type and ``Vermeidung von
Divisionen" (\cite{figismi}, \cite{giu-he92},
\cite{gihesa},
\cite{krick-pardo1},
\cite{krick-pardo:CRAS}; see also \cite{he-morg}).

Straight--line programs and arithmetic circuits  have their origins in the
design of semi--numerical algorithms and represent a link between numerical
analysis and computer science (see \cite{gathen:par},
\cite{joos:aaecc5},
\cite{Strassen:hand} and the references given there). In this section, we
slightly modify the notion of ``ordinary" arithmetic networks adding a special
type of computation nodes, called ``algebraic gates". These gates display
(generally algebraic) elements of ${\bar k}$ which are given as zeroes of a
precomputed univariate polynomial over $k$. We are now going to explain what we
mean by ``arithmetic network" and ``straight--line program with parameters in
${\bar k}$ (or $k$)". Let for the moment $K$ be any of the fields $k$ or ${\bar
k}$.

A piecewise rational function (over $K$) is a mapping $\varphi: {\bf
A}^n\longrightarrow {\bf A}^1$ such that there exists a partition of
${\bf A}^n$ in $K$--definable constructible subsets $\{C_j : 1\leq j
\leq M\}$, such that for any $1 \le j \le M$ there exists a rational function
$\phi_j \in K(X_1,\ldots,X_n)$ defined everywhere on $C_j$, which verifies~:

$$\varphi\mid_{_{\hbox {$C_j$}}}= {\phi_j}\mid_{_{\hbox {$C_j$}}}.$$ 

An ``{\it ordinary}" arithmetic network over $K$ is a device that evaluates
piecewise rational functions. However, in our applications, most of the
functions are piecewise polynomial, i.e. the
$\phi_j$ are polynomials belonging to $K[X_1,\ldots,X_n]$.

An arithmetic network $\Gamma$ over $K$ is a pair
$\Gamma=({\cal G}, Q)$, where ${\cal G}$ is directed acyclic graph, with
vertices (called nodes or gates of $\Gamma$) of indegree 0, 1, 2 or 3 and where
$Q$ is an assignment (labeling) of instructions and (piecewise) rational
functions (``intermediate results") to nodes (the labeling  will be specified
subsequently). The graph ${\cal G}$ contains $n$ nodes of indegree 0, labeled by
the variables $X_1,\ldots,X_n$, which are called the input nodes of $\Gamma$. 

We define the depth of a node $\nu$ of the graph ${\cal G}$ as the length of the
longest path joining $\nu$ with some input gate. Let us denote any node of
${\cal G}$ by a pair of integers $(i,j)$, where $i$ represents the depth of the
node and $j$ is an ordinal number assigned to the node which is given by some
ordering of the set of nodes of depth $i$ (see \cite{krick-pardo1}, 
\cite{Mo-Pardo:Net} for this type of encoding of arithmetic networks).

Let $(i,j)$ be a node of ${\cal G}$. If $(i,j)$ has indegree $0$ and does not
represent an input, the labeling $Q$ assigns to this node a constant from $K$
(which we call subsequently a {\it parameter} of
$\Gamma$). If $(i,j)$ has indegree
$2$ the labeling
$Q$ assigns to it an arithmetic operation of $k$. If $(i,j)$ has indegree $1$
the labeling
$Q$ assigns to it a {\it sign} (i.e. equal to zero) {\it test} and a {\it
boolean variable}. If
$(i,j)$ has indegree $3$ the labeling $Q$ assigns to it a {\it selector} which
makes a choice between two prefixed nodes of ${\cal G}$ according to the value
of a boolean variable associated with a third node. A node of
${\cal G}$ labeled in this manner is called a {\it gate} of
$\Gamma$. We have already introduced the input gates. A gate which is labeled by
a constant or an arithmetic operation of $K$ is called a {\it computation gate}
(in case of positive characteristic $p$ we consider also extraction of $p$--th
roots of elements of $K$ as arithmetic operations). The other gates are called
{\it sign} or {\it selector gates}, depending on their labeling. We suppose that
the labeling
$Q$ is ``meaningful", i.e. computation gates have ingoing edges (if any) coming
from computation or input gates, sign gates have an ingoing edge coming from a
computation gate , and selector gates have two ingoing edges coming from
computation or input gates and another ingoing edge coming from a sign gate.
Under this hypothesis the labeling $Q$ associates in a obvious way to each
computation gate
$(i,j)$ a rational function which we denote by
$Q_{i,j}$. Let
$(i,j)$ be a sign gate with ingoing edge coming from a node $(r,s)$ and let
$x$ be a point of ${\bf A}^n$ where $Q_{r,s}$ is defined. Then the boolean
variable
$B_{i,j}$ associated to $(i,j)$ takes the value $1$ or $0$ according to the
truth or falseness of the statement ``$Q_{i,j}(x)=0$". If
$(i,j)$ is a selector gate with two ingoing edges coming from computation gates
$(k,l)$ and
$(k',l')$ and if $x$ is a point of ${\bf A}^n$ where $Q_{k,l}$ and 
$Q_{k',l'}$ are defined and if the third ingoing edge of $(i,j)$ comes from the
sign gate
$(r,s)$ then we associate to
$(i,j)$ the field element $Q_{k,l}(x)$ or $Q_{k',l'}(x)$ according to the value
which takes in $x$ the boolean variable $B_{r,s}$ associated to the sign gate.
We consider the gates of $\Gamma$ with outdegree $0$ as {\it outputs} and we
suppose that  for each
$x \in {\bf A}^n$ there exists a ``stream" of ``consistent" paths (with all
rational functions defined in $x$) to a prefixed number of output gates.

A computation gate of $\Gamma$ is called {\it nonscalar} (with respect to $K$)
if the instruction associated to it corresponds to a multiplication of
nonconstant rational functions or to a division by a nonconstant rational
function. We assign to nonscalar computation gates unit costs whereas all other
gates are taken cost--free (in particular, computation gates of $\Gamma$
corresponding to
$K$--linear operations are free). In this way, we associate two complexity
measures to the arithmetic network $\Gamma$~:

\begin{itemize}

\item {\it nonscalar sequential time} or the {\it nonscalar size} of
      $\Gamma$, defined as the total number of nonscalar computation gates of
      $\Gamma$,

\item {\it nonscalar parallel time} or the {\it nonscalar depth} of
      $\Gamma$ defined as the longest oriented path of $\Gamma$ joining an
      input gate with an output gate when only nonscalar computation gates
      are taken into account for the ``length".  

\end{itemize} These two nonscalar complexity measures are fairly realistic with
respect to operation counting because $\Gamma$ can always be rearranged such
that the {\it total} number of gates and the {\it total} depth of $\Gamma$ are
bounded roughly by the square of its nonscalar size and depth respectively. The
intuitive meaning of the nonscalar size of $\Gamma$ is (sequential) running time
whereas the nonscalar depth is linked to rather mathematical quantities as
degree and height (in case 
$K\!:= \Q$) of the rational functions appearing as intermediate results of
$\Gamma$. The total depth of $\Gamma$ has a natural interpretation as minimal
storage space (see \cite{krick-pardo1} for more details on this subject).

A special case of arithmetic networks are those which have as single output a
boolean combination of sign gates. We call them  {\it decision networks}.
Another important case of special networks is represented by those which contain
only computation gates (no sign or selector gates). They are called {\it
straight--line programs} or {\it arithmetic circuits} and they compute (or
represent) rational functions belonging to
$K(X_1,\ldots,X_n)$. Often we call the (nonscalar) size of a straight--line
program its (nonscalar) {\it length}.

An arithmetic network or a straight--line program without any nonscalar division
(only divisions by non--zero elements of $K$ are allowed in this case) is called
{\it division free}. Division free straight--line programs compute polynomials of
$K[X_1,\ldots,X_n]$. Finally we say that a family of arithmetic networks with
nonscalar sequential time cost function $L$ is {\it well parallelizable} if its
nonscalar parallel time cost function $\ell$ behaves as $\ell=O(log_2 L)$. 

Since we are not considering other complexity measures than the nonscalar ones,
we simplify our terminology from now on saying just {\it size/sequential time}
or {\it depth/parallel time} with reference to the {\it nonscalar} cost measure.

Suppose now $K\!:=k$. We are going to extend the model of ``ordinary" arithmetic
networks (with parameters in $k$) introducing a new type of computation nodes,
called {\it algebraic gates}. These algebraic gates involve for given inputs
from $k$ elements of ${\bar k}$ which generally are algebraic over $k$ ( in
other words, these gates ``compute" algebraic functions in the inputs). In order
to explain the nature of these new gates let us assume that the arithmetic
network $\Gamma$ in question disposes over a second type of input nodes, labeled
by indeterminates, say $A_1,\ldots,A_m$ which are called {\it parameter
variables}. To understand this, remember that the concrete input for our
geometric problems (e.g. in Theorem \ref{afin}, \ref{tor} and
\ref{parametrizacion}) is always a division free straight--line program $\beta$
in
$k[X_1,\ldots,X_n]$ which represents the equations of the input system. The
parameter variables
$A_1,\ldots,A_m$ of $\Gamma$ are specialized for a concrete input circuit
$\beta$ into the values of the program parameters of
$\beta$, i.e. into values belonging to $k$. Therefore in all our applications
the number $m$ will be of order
$L^2$. Some computation nodes of the arithmetic circuit $\Gamma$ will therefore
depend exclusively on the parameter variables
$A_1,\ldots,A_m$. The labeling $Q$ assigns to these nodes a rational function of
$k(A_1,\ldots,A_m)$. We call such a node of $\Gamma$ {\it parameter gate}. Let
$T$ be a new indeterminate. An algebraic gate of
$\Gamma$ is now a node $(i,j)$ which has indegree and outdegree $N$, where $N$
is an arbitrary (but fixed) natural number. The ingoing edges of
$(i,j)$ are supposed to come from $N$ parameter gates of $\Gamma$, say
$(s_0,r_0),\ldots,(s_{N-1},r_{N-1})$, to which the labeling
$Q$ assigns polynomials $h_0,\ldots,h_{N-1} \in k[A_1,\ldots,A_m]$. Consider now
an arbitrary specialization $\alpha=(\alpha_1,\ldots,
\alpha_m) \in k^m$ of the parameter variables $A_1,\ldots,A_m$ in $k$ (such a
specialization represents in our applications always a concrete input circuit
$\beta$). The polynomials
$h_0,\ldots,h_{N-1}$ take in the argument $\alpha$ values
$\eta_0\!:=h_0(\alpha),\ldots,
\eta_{N-1}\!:=h_{N-1}(\alpha)$ which belong to $k$. The algebraic gate
$(i,j)$ assigns then to its $N$ outgoing edges the totality of the zeroes of the
{\it univariate} polynomial
$T^N+\eta_{N-1}T^{N-1}+\cdots +\eta_0 \in k[T]$ in arbitrary order (possibly
with repetitions if the polynomial is not separable). Let us observe that the
network
$\Gamma$ evaluates the polynomials $h_0,\ldots,h_{N-1}$ in the argument
$\alpha$ using possibly algebraic parameters which come from the previous use of
algebraic gates. Thus the network contains implicitly a straight--line program
representation of the values
$h_0(\alpha),\ldots,h_{N-1}(\alpha)$  which will be division free in all our
applications and which contains possibly algebraic parameters. If $\Gamma$ is an
arithmetic network which contains algebraic gates we shall say that $\Gamma$ has
{\it parameters in
${\bar k}$}.

\section{Compressing Straight--line Programs}
\label{compresion} In this section we show how we can compress in the algorithm
underlying Theorems
\ref{afin}, \ref{tor} and \ref{parametrizacion} some of the straight--line
programs which compute intermediate results. The principal outcome is the
following statement~:

\begin{lemma}[the affine case] 
\label{sinfor} Let be given polynomials $f_1,\ldots,f_i \in k[X_1,\ldots,X_n]$ 
and suppose that $f_1,\ldots,f_i$ are represented by a division free
straight--line program $\beta$ in
$k[X_1,\ldots,X_n]$ of length and depth $L$ and $\ell$ respectively. Assume that
the polynomials  $f_1,\ldots,f_i$ form a regular sequence in
$k[X_1,\ldots,X_n]$ and that they generate a radical ideal
$I\!:=(f_1,\ldots,f_i)$. Let $W\!:=V(f_1,\ldots,f_i)=V(I)$ be the affine variety
defined by $f_1,\ldots,f_i$ and denote by $\delta\!:= deg W$ the (geometric
affine) degree of $W$ and by
$r\!:= n-i$ its dimension (observe that by assumption the ideal $I$ is unmixed
and the variety
$W$ is equidimensional). Suppose that $X_1,\ldots,X_n$ are in Noether position
with respect to the variety $W$, the variables $X_1,\ldots,X_r$ being free and
assume that there is given by its coefficients a nonzero linear form
$u\in k[X_{r+1},\ldots,X_n]$ which represents a primitive element for $I$ (see
Section \ref{primitivo} below). The equations
$f_1,\ldots,f_i$, the variables $X_1,\ldots,X_n$ and the linear form $u$
determine uniquely the following mathematical objects~:  

\begin{itemize}

\item the minimal polynomial $q \in k[X_1,\ldots,X_r,u]$ of $u$ modulo
      the ideal $I$. This polynomial is monic and without loss of
      generality separable with respect to $u$ and satisfies $deg_u
      q=degq \le \delta$.  Let us call ${\bar \delta}:=deg_u q=deg q$. 

\item the (nonzero) discriminant $\rho \in k[X_1,\ldots,X_r]$ of $q$ and 
      polynomials $v_{r+1},\ldots,v_n \in k[X_1,\ldots,X_r,u]$ satisfying
      the conditions 
$$I_{\rho}=(q(u),\rho X_{r+1}- v_{r+1}(u),\ldots, \rho
      X_n- v_n(u))_{\rho}$$ $$max\{deg_u v_j\; ; \; r< j \le n\} <
      {\bar \delta} \le \delta.$$

\end{itemize}
 
Finally we assume that the polynomial $\rho$ and the coefficients of $q$ and
$v_{r+1},\ldots,v_n$ with respect to $u$ are given by a division free
straight--line program $\beta'$ in ${\bar k}[X_1,\ldots,X_r]$ having length and
depth $\Lambda$ and $\lambda$ respectively.

Under these assumptions there exists an arithmetic network with parameters in
${\bar k}$ which from the input circuits $\beta$ and $\beta'$ constructs a
division free straight--line program $\gamma$ in $k[X_1,\ldots,X_r]$ of length
$O((i^5+L) \delta^{11})$ and depth $O((log_2i+\ell) log_2
\delta)$ such that
$\gamma$ represents $\rho$ and  the coefficients of $q$ and
$v_{r+1},\ldots,v_n$  with respect to $u$. The size and depth of this network is
$O((i^5+L) \delta^{11})+\Lambda$ and $O((log_2i+\ell) log_2
\delta)+\lambda$ respectively.
\end{lemma}

For the proof of Lemma \ref{sinfor} we need the following fact~:
\medskip

\begin{remark}[the affine case]
\label{grado_uves} Let assumptions and notations be as in Lemma
\ref{sinfor}. Then the degrees of $\rho$ and of the  coefficients of
$v_{r+1},\ldots,v_n$ with respect to
$u$ are bounded by $2(deg q)^3 = 2 {\bar \delta}^3$, i.e. by $2 \delta^3$. 
\end{remark}

We show this remark later.

\medskip

{\it Proof of Lemma \ref{sinfor}.} 
 
Let us introduce the following notations~:

\begin{center}

$ B\!:= k[X_1,\ldots,X_n]/I  \; \; ; \;\;  A\!:= k[X_1,\ldots,X_r]$

\smallskip

$ K\!:=k(X_1,\ldots,X_r) \; \; ; \;\;  B'\!:=K[X_{r+1},\ldots,X_n]/(I), $

\smallskip
 
\end{center}
\noindent where $(I)$ denotes the ideal generated by the set $I$ (or by 
$f_1,\ldots,f_i$) in $K[X_{r+1},\ldots,X_n]$.

Recall that  the polynomials
$f_1,\ldots,f_i$ form a regular sequence in
$k[X_1,\ldots,X_n]$ and that the variables $X_1,\ldots,X_n$ are in Noether
position with respect to $W$, the variables $X_1,\ldots,X_r$ being free. Taking
this into account, we write 
$ \pi : W \longrightarrow {\bf A}^r$ for the finite and surjective morphism of
affine varieties induced by the coordinates $X_1,\ldots,X_r$.  Furthermore we
deduce from our assumptions (in particular from the assumption $q$ separable
with respect to $u$) that the finite dimensional 
$K$--algebra $B'$ is unramified and that $dim_K B'= deg_u q \le deg W=\delta$
holds (\cite{Iversen}). From this we infer that the jacobian

$$\Delta\!:= det \left( {{\partial f_k} \over {\partial X_j}}
\right)_{ 1
\le k \le i \atop r+1 \le j \le n}$$ is a nonzero divisor modulo $I$. Let us
observe that the polynomial $\Delta$ can be evaluated by a division free
straight--line program in $k[X_1,\ldots,X_n]$ of length 
$O(i^5+L)$ and depth $O(log_2 i +\ell)$. 
 
Let $\mu \in k[X_1,\ldots,X_r]$ be the constant term of the characteristic
polynomial of $\Delta$ modulo $I$. Since $\Delta$ represents a nonzero divisor
of $B$ we conclude $\mu \ne 0$. Furthermore we observe that $\mu$ can be
evaluated by a division free straight--line program in $k[X_1,\ldots,X_r]$ of
length 
$(i\delta L)^{O(1)}+\Lambda$ and depth $O(log_2(i\delta) + \ell)+\lambda$, and
so does the product $\rho \cdot \mu$. 

Using a correct test sequence (see \cite{heschnorr}) we are able to find in
sequential time $(i\delta L)^{O(1)}
\Lambda^2$ and parallel time $O(log_2(i \delta) + \ell)+\lambda$ a ``rational"
point $\eta=(\eta_1,\ldots,\eta_r)\in k^r$ which satisfies
$(\rho \cdot \mu) (\eta) \not= 0$ (thus we have $\rho(\eta) \not= 0$ and
$\mu(\eta) \not= 0$). From $\mu(\eta)\not= 0$ we deduce that the morphism
$\pi$ is unramified in the point $\eta$ what implies that
$\pi^{-1}(\eta)$ consists of exactly ${\bar \delta}$ nonsingular points  of $W$.
To be more precise, let $\pi^{-1}(\eta)= \{
\xi_1,\ldots,\xi_{{\bar \delta}} \} \subset W$ be the set of these points. Then
we infer from $\mu(\eta) \not= 0$ that $\Delta (\xi_l) \not= 0$ holds for any  $1
\le l \le {\bar \delta}$. 

>From $\pi(\xi_l)=\eta=(\eta_1,\ldots,\eta_r)$ we conclude that the point
$\xi_l$ has the form
$\xi_l=(\eta_1,\ldots,\eta_r,\xi_{r+1}^{(l)},\ldots,\xi_n^{(l)})
\in {\bar k}^n$ where the first $r$ coordinates $\eta_1,\ldots,\eta_r$ are
rational and independent of the index $l$ and the last $n-r$ coordinates
$\xi_{r+1}^{(l)},\ldots,\xi_n^{(l)}$ are algebraic (i.e. belong to ${\bar k}$
and not necessarily to $k$) and dependent on $l$. Let the linear form
$u$ have the form $u=\theta_{r+1}X_{r+1}+ \cdots+ \theta_nX_n$ with
$\theta_{r+1},\ldots,\theta_n \in k$. For $1 \le l \le {\bar \delta}$ we write
$\zeta_l\!:=u(\xi_l)=\theta_{r+1}\xi_{r+1}^{(l)}+ \cdots + \zeta_n \xi_n^{(l)}$
(observe that
$\zeta_l$ is typically an algebraic element of ${\bar k}$). Let us consider
$u$ as an indeterminate over $k$ and over $k[X_1,\ldots,X_r]$. In this sense
$q=q(X_1,\ldots,X_r,u)$ becomes a polynomial in $u$ with coefficients in
$k[X_1,\ldots,X_r]$ and $q(\eta,u)\!:=q(\eta_1,\ldots,\eta_r,u)$ becomes a
polynomial in $u$ with coefficients in $k$ (i.e. we have $q(\eta,u) \in k[u]$).
Since $q$ is monic in $u$ and of degree ${\bar \delta}$ we conclude that
$q(\eta,u)$ is monic and of degree ${\bar \delta}$ too. Furthermore $\rho(\eta)$
is the discriminant of $q(\eta,u)$ and from $\rho(\eta) \ne 0$ we deduce that
$q(\eta,u)$ has exactly ${\bar \delta}$ zeroes in ${\bar k}$ which are all
distinct. Let us analyze this fact more in detail~: for any polynomial 
$p \in k[X_1,\ldots,X_n]$ let us write
$p(\eta,X_{r+1},\ldots,X_n)\!:=p(\eta_1,\ldots,\eta_r,X_{r+1},\ldots,X_n)$ (thus
we have $p(\eta,X_{r+1},\ldots,X_n) \in k[X_{r+1},\ldots,X_n]$). We observe
that  $\{\xi_1\ldots,\xi_{{\bar \delta}} \}=\pi^{-1}(\eta)= \{\eta\}
\times V (f_1(\eta, X_{r+1},\ldots,X_n),\ldots,f_i(\eta, X_{r+1},\ldots,X_n))$
holds. From $I_{\rho}=(q(u),\rho X_{r+1}- v_{r+1}(u),\ldots, \rho X_n-
v_n(u))_{\rho}$ we deduce that for any $1 \le l
\le {\bar \delta}$ the identities
$$\xi_l=\left(\eta_1,\ldots,\eta_r,\displaystyle{v_{r+1}(\eta,\zeta_l)
\over
\rho(\eta)} ,\ldots,\displaystyle{v_n(\eta,\zeta_l) \over \rho(\eta)}
\right)$$ and $q(\eta,\zeta_l)=0$ hold (as before we write 
$v_j(\eta,u)\!:=v_j(\eta_1,\ldots,\eta_r,u) \in k[u]$ for $r < j \le n$). Since
the points $\xi_1,\ldots,\xi_{{\bar \delta}}$ are all different we conclude that
the values $\zeta_1,\ldots,\zeta_{{\bar \delta}} \in {\bar k}$ are distinct
zeroes of the polynomial $q(\eta,u) \in k[u]$. Since $q(\eta,u)$ is monic and of
degree
${\bar \delta}$ we conclude now that 
$$q(\eta,u)= \prod\limits_{l=1}^{{\bar \delta}}(u- \zeta_l)$$ holds. Thus the
values
$\zeta_1,\ldots,\zeta_{\delta}$ are exactly the zeroes of
$q(\eta,u)$.

Replacing in all polynomials $f_1,\ldots,f_i, \rho, q, v_{r+1},\ldots,v_n$ the
variables $X_1,\ldots,X_r$ by the new ones
$X_1-\eta_1,\ldots,X_r-\eta_r$ and leaving the remaining variables
$X_{r+1},\ldots,X_n$ and $u$ unchanged, we may assume without loss of generality
that $\eta=(0,\ldots,0) \in k^r$ holds. Let
$1\le l \le {\bar \delta}$. Recall that
$\xi_l=(\eta_1,\ldots,\eta_r,\xi_{r+1}^{(l)},\ldots,\xi_n^{(l)})=
(0,\ldots,0,\xi_{r+1}^{(l)},\ldots,\xi_n^{(l)})$, 
$$0=f_1(\xi_l)=f_1(\eta,\xi_{r+1}^{(l)},\ldots,\xi_n^{(l)})=
f_1(0,\ldots,0,\xi_{r+1}^{(l)},\ldots,\xi_n^{(l)})$$
$$\vdots$$ 
$$0=f_i(\xi_l)=f_i(\eta,\xi_{r+1}^{(l)},\ldots,\xi_n^{(l)})=
f_i(0,\ldots,0,\xi_{r+1}^{(l)},\ldots,\xi_n^{(l)})$$
 and 
$$0 \not=
\Delta(\xi_l)=
\Delta(\eta,\xi_{r+1}^{(l)},\ldots,\xi_n^{(l)}) =
\Delta(0,\ldots,0,\xi_{r+1}^{(l)},\ldots,\xi_n^{(l)})$$ holds. Thus the point 
$(\xi_{r+1}^{(l)},\ldots,\xi_n^{(l)}) \in {\bar k}^i$ is a nondegenerate zero of
the equation system given by the polynomials 
$f_1(\eta,X_{r+1},\ldots,X_n),\ldots,f_i(\eta,X_{r+1},\ldots,X_n) \in
k[X_{r+1},\ldots,X_n]$ which are in fact the polynomials
$f_1(0,\ldots,0,X_{r+1},\ldots,X_n),\ldots,f_i(0,\ldots,0,X_{r+1},\ldots,X_n)$.
>From Hensel's lemma (which represents a symbolic version of the Implicit
Function Theorem) we deduce that there exist formal power series
$R_{r+1}^{(l)},\ldots,R_n^{(l)} \in {\bar k}[[X_1,\ldots,X_r]]$ with 
$R_{r+1}^{(l)}(\eta)=\xi_{r+1}^{(l)},\ldots,R_n^{(l)}(\eta)=\xi_n^{(l)}$ such
that for
$R^{(l)}\!:=(X_1,\ldots,X_r,R_{r+1}^{(l)},\ldots,R_{n}^{(l)})$ the identities 
\begin{equation}
\label{fres0}   
        f_1(R^{(l)})=0,\ldots, f_i(R^{(l)})=0
\end{equation} hold in ${\bar k}[[X_1,\ldots,X_r]]$ (see \cite{Iversen} or
\cite{zar-sam}). Our next task is the construction of a ``sufficiently accurate
rational approximation" of the solution $R^{(l)}$ of the system given by the
polynomials
$f_1,\ldots,f_i$. This rational approximation will be represented by rational
functions ${\tilde R}_{r+1}^{(l)},\ldots,{\tilde R}_{n}^{(l)}\in {\bar
k}(X_1,\ldots,X_r)$ which are all defined in the point
$\eta=(0,\ldots,0)$ and can therefore be interpreted as power series belonging
to ${\bar k}[[X_1,\ldots,X_r]]$. For each
$r < j \le n$ the power series interpretation of the rational function
${\tilde R}_{j}^{(l)}$ will satisfy the congruence relation
\begin{equation}
\label{r-rtilde}  
    R_j^{(l)}-{\tilde R}_j^{(l)} \in (X_1,\ldots,X_r)^{2{\bar \delta}^3+1}
\end{equation} where $(X_1,\ldots,X_r)^{2{\bar \delta}^3+1}$ denotes the
($2{\bar \delta}^3+1$)--th power of the ideal generated by $X_1,\ldots$  \\
$\ldots,X_r$ in
${\bar k}[[X_1,\ldots,X_r]]$. The
$i$--tuple of rational functions $({\tilde R}_{r+1}^{(l)},\ldots,{\tilde
R}_{n}^{(l)})$ is defined by applying to the system $f_,\ldots,f_i$ at least
$3\lceil log_2{\bar \delta} \rceil +2$, i.e. roughly
$3\lceil log_2 \delta \rceil +2$ Newton iteration steps starting from the
particular nondegenerate solution
$\xi_l=(0,\ldots,0,\xi_{r+1}^{(l)},\ldots,\xi_n^{(l)})
\in {\bar k}^n$ (here we consider $f_1,\ldots,f_i$ as $i$--variate polynomials
depending on the variables $X_{r+1},\ldots,X_n$). In order to compute $({\tilde
R}_{r+1}^{(l)},\ldots,{\tilde R}_{n}^{(l)})$ we have to evaluate the polynomials
$f_k$ and their partial derivatives ${{\partial f_k} \over {\partial X_j}}$ for
$1 \le k \le i $ and $r< j \le n$ at $3\lceil log_2 \delta \rceil +2$ iteration
points. By \cite{baur-strassen} and \cite{krick-pardo1} (see also
\cite{morgenstern}) this can be done by a division free straight--line program
of size $O(Llog_2 \delta)$ and  depth $O(\ell log_2
\delta)$. Moreover we have to invert the jacobian matrix $\left( {{\partial f_k}
\over {\partial X_j}}
\right)_{ 1
\le k \le i \atop r < j \le n}$ at the same $3\lceil log_2 \delta \rceil+2$
iteration points. This costs additional $O(i^5 log_2 \delta)$ nonscalar
arithmetic operations organized in depth
$O(log_2 i log_2 \delta)$ (see \cite{berk}; observe also that these matrix
inversions require divisions). Thus $({\tilde R}_{r+1}^{(l)},\ldots,{\tilde
R}_{n}^{(l)})$ are represented by a straight--line program ${\tilde \gamma_l}$
in 
${\bar k}(X_1,\ldots,X_r)$ which contains as intermediate results only rational
functions which are defined in the point $\eta=(0,\ldots,0)$. The size and depth
of ${\tilde \gamma}_l$ are $O((i^5 + L) \log_2
\delta)$ and
$O((log_2 i +
\ell) \log_2 \delta)$ respectively.  

For each $1\le l \le {\bar \delta}$, let us write
${\tilde R}^{(l)}\!:=(X_1,\ldots,X_r,{\tilde R}_{r+1}^{(l)},\ldots,{\tilde
R}_{n}^{(l)})$. We are going to consider $u^{(l)}\!:=
u(R^{(l)})=\theta_{r+1}R_{r+1}^{(l)}+\cdots+\theta_nR_n^{(l)}$ and
${\tilde u}^{(l)}\!:=u({\tilde R}^{(l)})=\theta_{r+1}{\tilde
R}_{r+1}^{(l)}+\cdots+\theta_n{\tilde R}_{n}^{(l)}$. Observe that
$R_{r+1}^{(l)}(\eta)=\xi_{r+1}^{(l)},\ldots,R_n^{(l)}(\eta)=\xi_n^{(l)}$ implies
$u^{(l)}(\eta)=\zeta_l$. Without any nonscalar extra cost we may assume that
${\tilde \gamma}_l$ computes also the rational function 
${\tilde u}^{(l)} =\theta_{r+1}{\tilde R}_{r+1}^{(l)}+\cdots+\theta_n{\tilde
R}_{n}^{(l)}$. Let ${\tilde \gamma}$ be the straight--line program in ${\bar
k}(X_1,\ldots,X_r)$ which we obtain by joining all circuits
${\tilde \gamma}_1,\ldots,{\tilde \gamma}_{{\bar \delta}}$. The straight--line
program
${\tilde \gamma}$ computes ${\tilde u}^{(1)},\ldots,{\tilde u}^{({\bar
\delta})}$ and has  size 
$O((i^5 + L) \delta \log_2 \delta)$ and depth $O((log_2 i +
\ell) \log_2 \delta)$. Observe that all the intermediate results of
${\tilde
\gamma}$ are defined in the point $\eta=(0,\ldots,0)$.

Let $1 \le l \le {\bar \delta}$. Taking into account that
$q$ is the minimal polynomial of $u$ modulo the ideal $I=(f_1,\ldots,f_i)$ we
deduce from (\ref{fres0}) and the assumptions of the lemma the identity 
$q(X_1,\ldots,X_r,u^{(l)})=0$. Moreover (\ref{r-rtilde}) implies that the
congruence relation
\begin{equation}
\label{u-utilde}
    u^{(l)}-{\tilde u}^{(l)} \in (X_1,\ldots,X_r)^{2{\bar \delta}^3+1}
\end{equation} holds in the power series ring ${\bar k}[[X_1,\ldots,X_r]]$.

Recall that we consider $q$ as a polynomial in $u$ with coefficients in
$k[X_1,\ldots,X_r]$. Moreover $q$ is monic in $u$ and of total degree
${\bar \delta}$. Thus the polynomial $q$ has the form $q=\sum\limits_{0\le m
\le {{\bar \delta}}} q_{m}u^m$ with $q_m \in k[X_1,\ldots,X_r]$ and $deg q_m \le
{\bar \delta}$ for $ 0 \le m \le {\bar \delta}$. Furthermore we have
$q_{_{\bar \delta}}=1$. On the other hand we know already that
$q(X_1,\ldots,X_r,u^{(l)})=0$ holds for any $1
\le l \le {\bar \delta}$.  From $u^{(l)}(\eta)=\zeta_l$ we deduce that the power
series
$u^{(1)},\ldots,u^{({\bar \delta})}$ are all distinct. Since $q$ is monic and of
degree
${\bar \delta}$ in $u$ and  $u^{(1)},\ldots,u^{({\bar \delta})}$  are distinct
zeroes of
$q$ we conclude that in ${\bar k}[[X_1,\cdots,X_r]][u]$ the identity  
\begin{equation}
\label{qyraices} 
      q=\prod\limits_{1\le l \le {\bar \delta}} (u-u^{(l)})
\end{equation} holds.  For $0 \le m \le {\bar \delta}$ let us denote by
$\sigma_m$ the
$m$--th elementary symmetric function in ${\bar \delta}$ arguments. From
(\ref{qyraices}) we deduce $q_m=(-1)^{{\bar \delta} - m}
\sigma_m(u^{(1)},\ldots,u^{({\bar \delta})})$. In combination with
(\ref{u-utilde}) this implies that for $ 0 \le m < {\bar \delta}$ the congruence
relation
\begin{equation}
\label{pm-sigmam}
  q_m - (-1)^{{\bar \delta} - m} \sigma_m({\tilde u}^{(1)},\ldots,{\tilde
u}^{({\bar \delta})}) \in (X_1,\ldots,X_r)^{{\bar \delta}+1}
\end{equation}  holds in the power series ring ${\bar k}[[X_1,\ldots,X_r]]$. As
for $ 0 \le m\le {\bar \delta}$ the polynomial $q_m \in k[X_1,\ldots,X_r]$ has
degree at most ${\bar \delta}$ we deduce from (\ref{pm-sigmam}) that the power
series expansion of 
$(-1)^{{\bar \delta} - m} \sigma_m({\tilde u}^{(1)},\ldots,{\tilde u}^{({\bar
\delta})})
\in {\bar k} [[X_1,\ldots,X_r]]$ coincides up to degree ${\bar \delta}$ with
$q_m$. Combining the straight--line program ${\tilde \gamma}$ which computes the
rational functions
${\tilde u}^{(1)},\ldots,{\tilde u}^{({\bar \delta})}$ with a fast and
well--parallelizable algorithm for the evaluation of the set of elementary
symmetric functions $\{ \sigma_m\; ; \; 0 \le m \le {\bar \delta} \}$ (see
\cite{strassen2}, Satz 3.1. and \cite{jaja}, Chapter 8, Exercise 8.15) we obtain
an arithmetic circuit $\gamma_0$ in ${\bar k}(X_1,\ldots,X_r)$ which for $0 \le
m < {\bar \delta}$ computes all rational functions $(-1)^{{\bar \delta} - m}
\sigma_m({\tilde u}^{(1)},\ldots,{\tilde u}^{({\bar \delta})})$. The circuit
$\gamma_0$ has  size $O((i^5 + L)
\delta  \log_2  \delta)$ and depth $O((log_2 i + \ell) \log_2 \delta)$ and all
its intermediate results are rational functions of ${\bar k}(X_1,\ldots,X_r)$
which are defined in $\eta=(0,\ldots,0)$. Taking this last observation into
account we apply to the circuit $\gamma_0$ the well--parallelizable Vermeidung
von Divisionen technique contained in the proof of
\cite{krick-pardo1}, Proposition 21. (With respect to parallelism we observe
here that only one division by a suitable power of $\Delta$ is necessary at the
very end of the procedure. If one is interested only in the sequential aspect of
this technique one may apply directly the simpler algorithm underlying 
\cite{stras:verm}, Satz 2). In this way we obtain a {\it division free}
straight--line program
$\gamma_1$ in
${\bar k}[X_1,\ldots,X_r]$ which for $0 \le m < {\bar \delta}$ computes the power
series expansion in ${\bar k}[[X_1,\ldots,X_r]]$ of all rational functions $
(-1)^{{\bar \delta} - m}
\sigma_m({\tilde u}^{(1)},\ldots,{\tilde u}^{({\bar \delta})})$ {\it up to degree
${\bar \delta}$}. Since  for $ 0 \le m < {\bar \delta}$ these truncated power
series expansions coincide with the polynomials $q_m \in k[X_1,\ldots,X_r]$ we
conclude that $\gamma_1$ represents these polynomials and hence the coefficients
of $q$ in the representation  $q=\sum\limits_{0\le m \le {\bar \delta}}
q_{m}u^m$. The  size and depth of $\gamma_1$ are $O((i^5 + L)
\delta^3 \log_2
\delta)$ and $O((log_2 i + \ell) \log_2 \delta)$ respectively.

Once given the straight--line program representation $\gamma_1$ of the
coefficients $q_m$  of the polynomial $q$ we are able to compute the
discriminant $\rho$ using additional $O(\delta^5)$  arithmetic operations
organized in depth $O(log_2
\delta)$. Therefore we can extend the circuit $\gamma_1$ to a division free
straight--line program $\gamma_2$ in ${\bar k}[X_1,\ldots,X_r]$ of roughly the
same  size and depth as $\gamma_1$, such that
$\gamma_2$ computes also the discriminant $\rho$.

In order to finish the proof of Lemma \ref{sinfor} we have to find a
straight--line program representation of the coefficients of the polynomials 
$v_{r+1},\ldots,v_n$ with respect to the variable $u$ such that the resulting
circuit has about the same  size and depth as
$\gamma_2$. By the way we shall assume without loss of generality that the
circuit $\gamma_1$ (and hence
$\gamma_2$) contains the circuit ${\tilde
\gamma}$.

Let $ r< j \le n$ and let $v_j=\sum\limits_{0\le m < {\bar \delta}}
a_{m}^{(j)}u^m$ with $a_{m}^{(j)} \in k [X_1,\ldots,X_r]$.From (\ref{fres0}) and
the assumptions of the lemma we deduce that for any  $1
\le l \le {\bar \delta}$ the identity
\begin{equation}
\label{seis}  
       v_j(u^{(l)}) -\rho R_j^{(l)}=0 
\end{equation} holds. Taking into account the degree bound $deg_u v_j < {\bar
\delta}$ we interpret the identities (\ref{seis}) as an inhomogeneous
${\bar \delta} \times {\bar \delta}$ linear equation system for the coefficient
vector
$(a_{0}^{(j)},\ldots,a_{{\bar \delta} -1}^{(j)})$ of $v_j$. The matrix of this
equation system is the ${\bar \delta} \times {\bar \delta}$ Vandermonde matrix
corresponding to the
${\bar \delta}$ (distinct) power series $(u^{(1)},\ldots,u^{({\bar \delta})})$
and its inhomogeneous part is the row vector we get  transposing the
${\bar \delta}$--tuple
$(\rho R_j^{(1)},\ldots,\rho R_j^{({\bar \delta})})$. From (\ref{r-rtilde}) and
(\ref{u-utilde}) and the assumptions of the lemma we deduce as before that for
any $ 1 \le l
\le {\bar \delta}$ the congruence relation 
\begin{equation}
\label{siete} v_j({\tilde u}^{(l)})-\rho {\tilde R}_j^{(l)} \in
(X_1,\ldots,X_r)^{2{\bar \delta}^3+1}
\end{equation}  holds in  ${\bar k}[[X_1,\ldots,X_r]]$. We interpret now the
congruence relations (\ref{siete}) as an inhomogeneous ${\bar \delta} \times
{\bar \delta}$ linear equation system whose (unique) solution is a rational
approximation to the coefficient vector $(a_{0}^{(j)},\ldots,a_{{\bar \delta}
-1}^{(j)})$ of $v_j$. More precisely we extend the circuit $\gamma_2$ which
computes $\rho$ and for
$ 1 \le l \le {\bar \delta}$ the rational functions ${\tilde u}^{(l)}$ and
${\tilde R}_j^{(l)}$ to a straight--line program $\gamma'_j$ in ${\bar
k}(X_1,\ldots,X_r)$ as follows~: for $ 1 \le l \le {\bar \delta}$ the circuit
$\gamma'_j$ contains as intermediate results the rational functions $\rho
{\tilde R}_j^{(l)}$. Moreover
$\gamma'_j$ computes rational functions ${\tilde a}_{0}^{(j)},\ldots,{\tilde
a}_{{\bar \delta} -1}^{(j)}$ which are obtained by multiplying the ${\bar
\delta}$--tuple
$(\rho R_j^{(1)},\ldots,\rho R_j^{({\bar \delta})})$ by the inverse of the
Vandermonde matrix corresponding to the (distinct) rational functions
$({\tilde u}^{(1)},\ldots,{\tilde u}^{({\bar \delta})})$. We observe that all
intermediate results of $\gamma'_j$ are rational functions of ${\bar
k}(X_1,\ldots,X_r)$ which are defined in the point $\eta=(0,\ldots,0)$. Moreover
$\gamma'_j$ has  size and depth
$O((i^2+L)\delta^3log^\delta + \delta^5)$ and $O((log_2i + \ell) log_2
\delta)$ respectively. From (\ref{seis}) and (\ref{siete}) we deduce that for
any $0 \le m < {\bar \delta}$ the congruence relation 

\begin{equation}
\label{ocho}
     a_m^{(j)}-{\tilde a}_m^{(j)} \in (X_1,\ldots,X_r)^{2{\bar \delta}^3+1}
\end{equation} holds in ${\bar k}[[X_1,\ldots,X_r]]$.

By Remark \ref{grado_uves} the degree of any polynomial $a_m^{(j)} \in
k[X_1,\ldots,X_r]$ is bounded by $2{\bar \delta}^3$. Therefore (\ref{ocho})
implies that the power series expansion of the rational function ${\tilde
a}_m^{(j)}$ in
${\bar k}[[X_1,\ldots,X_r]]$ coincides up to degree $2{\bar \delta}^3$ with $
a_m^{(j)}$. We join now all the circuits
$\gamma'_{r+1},\ldots,\gamma'_n$ and apply to them the same Vermeidung von
Divisionen technique as before. In this way we obtain a division--free
straight--line program $\gamma$ in ${\bar k}[X_1,\ldots,X_r]$ which computes for
$r < j \le n$ and $0 \le m < {\bar \delta}$ the power series expansions in ${\bar
k}[[X_1,\ldots,X_r]]$ of all rational functions
${\tilde a}_m^{(j)}$ {\it up to degree $2{\bar \delta}^3$}. Since for $r < j \le
n$ and
$0 \le m < {\bar \delta}$ this truncated power series expansions coincide with
the polynomials
$a_m^{(j)}$, the circuit $\gamma$ represents all coefficients with respect to
$u$ of the polynomials $v_{r+1},\ldots,v_n$. Without loss of generality  we may
assume that the straight--line program $\gamma$ extends the circuit
$\gamma_2$. Thus finally the straight--line program $\gamma$ computes
$\rho $ and all the coefficients of $q$ and $v_{r+1},\ldots,v_n$. with respect
to $u$. The  size and depth of $\gamma$ are $ O((i^5 +L)
\delta^{9} log_2 
\delta + \delta^{11})= O((i^5+L) \delta^{11})$ and $O((log_2i + \ell) log_2
\delta)$ respectively. 
\cqfd

\medskip

We shall make use of Lemma \ref{sinfor} in the proof of Theorem
\ref{afin} and Theorem \ref{parametrizacion} $(i)$ which deal with the 
zerodimensional algorithmic elimination problem in the {\it affine case}. For
the same problem in the {\it toric case}, i.e. for Theorem \ref{tor} and Theorem
\ref{parametrizacion} $(ii)$, we need a slightly different version of Lemma
\ref{sinfor}.

This is the content of the next statement.

\begin{lemma}[the toric case] 
\label{sinfortorico}Let be given polynomials $f_1,\ldots,f_i \in
k[X_1,\ldots,X_n]$  and suppose that $f_1,\ldots,f_i$ are represented by a
division free straight--line program $\beta$ in
$k[X_1,\ldots,X_n]$ of length and depth $L$ and $\ell$ respectively. Let also
$g:= {\prod\limits_{j=1}^n X_j}$. Assume that the polynomials  $f_1,\ldots,f_i$
form a toric complete intersection ideal
$I\!:=(f_1,\ldots,f_i)$ whose localization $I_g$ is radical. Let
$W^{*}$ be the union of the toric irreducible components of the affine variety
$V(f_1,\ldots,f_i)=V(I)$ defined by
$f_1,\ldots,f_i$ and denote by
$\delta^{*}\!:= deg W^{*}=deg^{*} V(f_1,\ldots,f_i)$ the (geometric) degree of
$W^{*}$ and by
$r\!:= n-i$ its dimension (observe that by assumption the localized ideal
$I_g$ is unmixed and the ``toric" variety
$W^{*}$ is equidimensional). We assume that $W^{*}$ is not empty, i.e.
$\delta^{*} > 0$. Suppose that
$X_1,\ldots,X_n$ are in Noether position with respect to the variety
$W^{*}$, the variables
$X_1,\ldots,X_r$ being free and assume that there is given by its coefficients a
nonzero linear form
$u\in k[X_{r+1},\ldots,X_n]$ which represents a primitive element for
$I_g$ (see Section \ref{primitivo} below). The equations
$f_1,\ldots,f_i$, the variables $X_1,\ldots,X_n$ and the linear form $u$
determine uniquely the following mathematical objects~:  

\begin{itemize}

\item a polynomial $q^{*} \in k[X_1,\ldots,X_r,u]$ which is monic
      and without loss of generality separable with respect to $u$ and satisfies
      $deg_u q^{*}=degq^{*} \le \delta^{*}$  

\item the (nonzero) discriminant $\rho^{*} \in k[X_1,\ldots,X_r]$ of
      $q^{*}$ and polynomials $v^{*}_{r+1},\ldots,v^{*}_n \in
      k[X_1,\ldots,X_r,u]$ satisfying
      the conditions 
$$I_{\rho^{*} g} =(q^{*}(u),\rho^{*}
      X_{r+1}-
      v^{*}_{r+1}(u),\ldots, \rho^{*} X_n- v^{*}_n(u))_{\rho^{*} g}$$ 
$$max\{deg_u v^{*}_j\; ; \; r< j \le
      n\} < deg q^{*} \le \delta^{*}.$$

\end{itemize}
 
Finally we assume that the polynomial $\rho^{*}$ and the coefficients of
$q^{*}$ and
$v^{*}_{r+1},\ldots,v^{*}_n$ with respect to $u$ are given by a division free
straight--line program $\beta^{*}$ in ${\bar k}[X_1,\ldots,X_r]$ having length
and depth $\Lambda^{*}$ and $\lambda^{*}$ respectively.

Under these assumptions there exists an arithmetic network with parameters in
${\bar k}$ which from the input circuits $\beta$ and $\beta^{*}$ constructs a
division free straight--line program $\gamma^{*}$ in
${\bar k}[X_1,\ldots,X_r]$ of length
$O((i^5+L) {\delta^{*}}^{11})$ and depth $O((log_2i+\ell) log_2
\delta^{*})$ such that
$\gamma^{*}$ represents $\rho^{*}$ and  the coefficients of $q^{*}$ and
$v^{*}_{r+1},\ldots,v^{*}_n$  with respect to $u$. The size and depth of this
network is $O((i^5+L) {\delta^{*}}^{11})+\Lambda^{*}$ and
$O((log_2i+\ell) log_2
\delta^{*})+\lambda^{*}$ respectively.

\noindent This statement remains true if ``toricity" is understood with respect
to a linear change of the original variables $X_1,\ldots,X_n$ and the polynomial
$g=\prod\limits_{j=1}^n X_j$ is replaced by the product of the new variables.
\end{lemma}

Since the proof of Lemma \ref{sinfor} is completely local, only minor changes at
its very beginning and the application of a suitable toric  version of Remark
\ref{grado_uves} at its very end are necessary to convert it into a proof of
Lemma
\ref{sinfortorico}. The essential steps and arguments are almost textually the
same. The most subtle point is the {\it toric} version of Remark
\ref{grado_uves} we need at the end of the proof. For the statement of this
toric version we replace just in Remark \ref{grado_uves} the polynomials 
$\rho$, $v_{r+1},\ldots,v_n$ by $\rho^{*}$, 
$v^{*}_{r+1},\ldots,v^{*}_n$ and the degree bound by 
$2 (deg q^{*})^3$, i.e. by $ 2 (\delta^{*})^3$. We shall come back to the proof
of that at the end of this section. In order to avoid repetitive argumentation we
omit the proof of Lemma
\ref{sinfortorico}. 

The following Remark \ref{fact_gates} refers to work in progress. It  answers
the question whether computing with algebraic elements of
${\bar k}$ is really necessary in the algorithm underlying Lemma
\ref{sinfor} and 
\ref{sinfortorico}. On the other hand one might ask whether to this algorithm
corresponds a counterpart in the bit model of boolean circuits in case
$k\!:=\Q$. We postpone the answer to this second question to the end of the
paper (see Section \ref{conclusiones}).

Suppose now that $k$ is a hilbertian field for which a factorization algorithm
for univariate polynomials is available at {\it ``moderate arithmetic costs"}.
This means that we may add (in a computationally reasonable way) to our
arithmetic networks over $k$ special {\it factorization gates} for univariate
polynomials of degree $D$, where
$D$ is an arbitrary (but fixed) natural number. Our new complexity model takes
such a gate into account at costs of
$D^{O(1)}$ with respect to sequential time (this is quite realistic in view of
\cite{LLL}) and at costs of
$O(log_2D)$ with respect to parallel time (this is rather cheap). We call such a
base field $k$ {\it hilbertian with univariate polynomial factorization at
moderate costs} and the corresponding algebraic complexity model is called {\it
arithmetic network (over $k$) with factorization gates}. Observe that this model
is {\it nonuniform}  by the way we have introduced it. However in view of
\cite{ha-mon} and the evidence that ``efficient" (rather than ``effective")
versions of Hilbert's Irreducibility Theorem must exist (see e.g. \cite{debes},
\cite{fried}, \cite{sprindzak} and
\cite{yasumoto}), we may hope to obtain at least a reasonable randomization (if
not uniformization) of the algorithm underlying the next statement which we
pronounce without proof.

\begin{remark}
\label{fact_gates} Suppose that $k$ is a hilbertian field with univariate
factorization at moderate costs. Suppose that there is given the same situation 
(with the same notations) as in Lemma \ref{sinfor} and
\ref{sinfortorico} and make the additional assumption that the circuits
$\beta'$ and $\beta^{*}$ use only parameters from the base field $k$. Then,
subject to the following modifications, the same conclusions as in Lemma
\ref{sinfor} and \ref{sinfortorico}
 are true~:
\begin{itemize}
\item the output network in question is an arithmetic network over $k$
      with factorization gates
\item the straight--line programs $\gamma$ and $\gamma^{*}$ are division free
circuits in $k[X_1,\ldots,X_r]$ (which use only parameters from $k$).
\end{itemize}
\end{remark}

We finish this section showing the assertion of Remark \ref{grado_uves} and of
its toric counterpart.

\medskip

{\it Proof of Remark \ref{grado_uves}}. Let assumptions and notations be as in
Lemma \ref{sinfor} and its proof. We shall consider the $k$--algebras $C\!:=
k[X_1,\ldots,X_r,u]/q(u)$ and

$$C_{\rho}\!:=k[X_1,\ldots,X_r]_{\rho}[u]/(q(u))_{\rho}$$ and compare them to
the $k$--algebra $B\!:=k[X_1,\ldots,X_n]/(f_1,\ldots,f_i)$ and its localization
$$B_{\rho}\!:=k[X_1,\ldots,X_r]_{\rho}[X_{r+1},\ldots,X_n]/
(f_1,\ldots,f_i)_{\rho}.$$ Let us denote by $\bar{}$ the canonical $k$--algebra
homomorphism which associates to each element in $k[X_1,\ldots,X_n]$,

\noindent $k[X_1,\ldots,X_r]_{\rho}[X_{r+1},\ldots,X_n]$, $k[X_1,\ldots,X_r,u]$
or
$k[X_1,\ldots,X_r]_{\rho}[u]$ its residue class in $B$, $B_{\rho}$,  $C$  or
$C_{\rho}$. Observe that the condition
$(f_1,\ldots,f_i)_{\rho}=I_{\rho}=(q(u),\rho X_{r+1}-v_{r+1}(u),\ldots,
\rho X_n-v_n(u))_{\rho}$ induces a $k$--algebra isomorphism $\tau:B_{\rho}
\longrightarrow C_{\rho}$ which leaves $k[X_1,\ldots,X_r]_{\rho}$ fixed and
which satisfies the condition $\tau({\bar X_{r+1}})=
\displaystyle{v_{r+1}({\bar u}) \over
\rho},\ldots,\tau({\bar X_n})= \displaystyle{v_n({\bar u}) \over
\rho}$. Since the polynomials $f_1,\ldots,f_i$ form a regular sequence in
$k[X_1,\ldots,X_n]$ and the variables $X_1,\ldots,X_n$ are in Noether position
with respect to the ideal $(f_1,\ldots,f_i)$, the $k$--algebra $B$ is a free
$k[X_1,\ldots,X_r]$--module of rank ${\bar \delta}$ (see e.g.
\cite{gihesa}, Lemma 3.3.1). The same is true for the $k$--algebra $C$. Hence
$B_{\rho}$ and $C_{\rho}$ are free
$k[X_1,\ldots,X_r]_{\rho}$--modules of rank ${\bar \delta}$. The polynomial is
$q$ by assumption separable with respect to the variable $u$. Thus we conclude
that $B_{\rho}$ and $C_{\rho}$ are unramified
$k[X_1,\ldots,X_r]_{\rho}$--algebras (this fact has already been used in the
proof of Lemma \ref{sinfor}). Therefore the (ordinary) traces which map $B$ and
$C$ onto $k[X_1,\ldots,X_r]$ are nonzero and the corresponding traces of
$B_{\rho}$ and $C_{\rho}$ are nondegenerate (observe that $\rho$ is the
discriminant of the
$k[X_1,\ldots,X_r]$--module basis of $C$ given by the elements ${\bar 1},{\bar
u},\ldots,{\bar u}^{{\bar \delta}-1}$). Since the $k$--algebra isomorphism
$\tau$ is
$k[X_1,\ldots,X_r]_{\rho}$--linear it leaves the traces fixed. We denote them
therefore by the same symbol {\it ``Tr"} (see
\cite{kunz}, Appendix F as reference for traces in noetherian algebras over
rings). Let $r < j \le n$. As in the proof of Lemma \ref{sinfor} let us write
$v_j=\sum\limits_{0\le m < {\bar \delta}} a_{m}^{(j)}u^m$ with
$a_{m}^{(j)} \in k[X_1,\ldots,X_r]$. Observe that $\rho {\bar X_j}=v_j({\bar u})
=\sum\limits_{0\le m < {\bar \delta}} a_{m}^{(j)}{\bar u}^m$ holds in $B$. For
$0 \le m' < {\bar \delta}$ consider the polynomial 
\begin{equation}
\label{n9} c_{m'}^{(j)}\!:= \rho Tr({\bar X_j}{\bar u}^{m'})= \sum\limits_{0\le
m <
  {\bar \delta}} a_{m}^{(j)} Tr({\bar u}^{m+m'})
\end{equation} which belongs to $k[X_1,\ldots,X_r]$. 

Observe that $Tr({\bar X_j}{\bar u}^{m'})$ is also a polynomial of
$k[X_1,\ldots,X_r]$ which satisfies by \cite{saso}, Theorem 13 the degree bound
$deg Tr({\bar X_j}{\bar u}^{m'}) \le deg W deg(X_j u^{m'})= (m'+1) {\bar \delta}
\le {\bar \delta}^2$. Thus $c_{m'}^{(j)}$ is divisible in
$k[X_1,\ldots,X_r]$ by $\rho$ and $\displaystyle{c_{m'}^{(j)} \over
\rho}$ satisfies the degree bound 
 \begin{equation}
\label{n10}
  deg \displaystyle{c_{m'}^{(j)} \over \rho} \le {\bar \delta}^2.
\end{equation}

In the same way we see that for any $0 \le m < {\bar \delta}$ the trace
$Tr({\bar u}^{m+m'})$ is an element of $k[X_1,\ldots,X_r]$ which satisfies the
degree condition

 \begin{equation}
\label{n11}
 deg Tr({\bar u}^{m+m'}) \le 2{\bar \delta} ({\bar \delta} -1).
\end{equation}

For $j$ fixed and $m$ verifying $0 \le m < {\bar \delta}$ we interpret the
equations (\ref{n9}) as an inhomogeneous $ {\bar \delta} \times {\bar \delta}$
linear equation system for the coefficient vector
$(a_0^{(j)},\ldots,a_{{\bar \delta}-1}^{(j)})$. The matrix of this equation
system is
$\left( Tr({\bar u}^{m+m'}) \right) _{0 \le m,m' < {\bar \delta}}$ whose entries
belong to $k[X_1,\ldots,X_r]$. The determinant of the equation system is
therefore
$\rho$, which is the discriminant of the $k[X_1,\ldots,X_r]$--module basis
$1,{\bar u},\ldots,{\bar u}^{{\bar \delta} -1}$ of $C$. The inhomogeneous vector
of the equation system is obtained transposing the ${\bar \delta}$--tuple
$(c_0^{(j)},\ldots,c_{{\bar \delta} - 1}^{(j)})$.

>From Cramer's rule and the degree bounds (\ref{n10}) and (\ref{n11}) we conclude
therefore 
$$deg a_m^{(j)} \le 2 {\bar \delta}^3$$
 for any $r < j \le n$ and $0 \le m < {\bar \delta}$. Similarly we deduce from
(\ref{n11}) the degree bound $deg \rho \le 2{\bar \delta}^3 \le 2\delta^3$.
\cqfd

\medskip
 
The toric version of Remark \ref{grado_uves} can be proved essentially by the
same method. One has only to replace the $k$--algebra $B$ by the (reduced)
coordinate ring of the ``toric" variety $W^{*}$. Of course this coordinate ring
is not a free $k[X_1,\ldots,X_r]$--module anymore, but it is still a finite
faithful module over  $k[X_1,\ldots,X_r]$, i.e. an integral extension of
$k[X_1,\ldots,X_r]$. This is sufficient for traces to be well defined and for
\cite{saso}, Theorem 13 to be true also in this new context. The rest of the
proof of the toric version of Remark \ref{grado_uves} is textually the same as
in the affine case.

Let us also observe that Remark \ref{grado_uves} can be formulated with a
slightly better degree bound of ${\bar \delta}^2+1$, i.e. of $\delta^2+1$in the
affine case and
$(degq^{*})^2+1$, i.e of $\delta^{{*}2}+1$ in the toric case. This can be easily
deduced from our elementary arguments in Section \ref{alg} taking into account
that all the polynomials involved are homogeneous.

\section{Technical Lemmas}
\label{lemas}
\subsection{Squarefree Representation and Greatest Common Divisor Computation of
Univariate Polynomials}
\label{subsecuno} Let $R$ be an integral and factorial $k$--algebra with
fraction field $K$, $T$ a new variable,
$D$ a fixed natural number and $P,Q \in R[T]$ two polynomials of formal degree
$D$. We think that $P$ and $Q$ are given by their coefficient vectors each of
length $D$. In the statements which follow now we refer to (ordinary) arithmetic
networks and straight--line programs over $R$ which receive the coefficient
vectors of $P$ and $Q$ as inputs. Furthermore greatest common divisors (gcd's)
will always be taken with respect to the principal ideal domain $K[T]$ although
they are represented in $R[T]$ (thus they are not unique up to units in R). The
proofs of the following well-known lemmas can be found in \cite{krick-pardo1}.

\begin{lemma}
\label{gcd} There exists a division free well parallelizable arithmetic network
$\Gamma$ of size  $O(D^{6})$ which from the coefficients of
$P$ and $Q$ computes the coefficients of a greatest common divisor (belonging to
$R[T]$) of these polynomials. 
\end{lemma}

The greatest common divisor computed by the algorithm underlying this Lemma is a
(not necessarily primitive) polynomial of $R[T]$ which we denote by
$gcd(P,Q)$. A similar well parallelizable ($O(s^5D^6)$ sequential time )
complexity bound holds for the computation of the simultaneous greatest common
divisor of more than two (say
$s$) polynomials of degree at most
$D$ belonging to
$R[T]$.

Sometimes the second polynomial $Q$ is given by a division free straight--line
program in $R[T]$ of length $L$ and depth $\ell$ instead by its coefficients. In
this case we have

\begin{lemma}
\label{gcd2}
 There exists a division free arithmetic network of size $LD^{O(1)}$ and depth
$O(log_2(D)+ \ell)$ which computes the coefficients of the greatest common
divisor of $P$ and $Q$. 
\end{lemma} 

The algorithm underlying Lemma \ref{gcd2} uses linear algebra subroutines which
deal with square matrices over $R$ of size at most
$2D-1$. 

As before this lemma can be generalized to the case of the computation of the
simultaneous greatest common divisor of more than two (say $s$) polynomials, one
of them given by its coefficient vector and having degree at most $D$ and the
others given by a division free straight--line program in $R[T]$ of length
$L$ and depth $\ell$. The outcome is then a division free arithmetic network of
size
$(L+s)D^{O(1)}$ and depth $O(log_2(D)+ \ell)$ which computes the coefficients of
the greatest common divisor in question. Again the linear algebra subroutines
deal only with square matrices of size at most $2D-1$. Lemma
\ref{gcd} can be used in order to compute a separable (and hence in $K[T]$
squarefree) polynomial
${\tilde P} \in R[T]$ which has the same zeroes as $P$ (we call such a
polynomial ${\tilde P}$ a {\it separable representation} of
$P$). If the characteristic of
$k$ (and hence of $R$ and $K$) is zero this is a immediate consequence of the
next lemma putting ${\tilde P}\!:=P^{*}$. In case of positive characteristic we
need a more refined analysis of the situation.  

Observe that the coefficients of the derivative $P'$ of $P$ are immediately
obtained from the (given) coefficients of $P$. Therefore we have

\begin{lemma} \label{sqfree} Suppose that the derivative $P'$ is nonzero. Then
there exists a division free well parallelizable arithmetic network of size
$O(D^{6})$ which from the coefficients of $P$ computes the coefficients of a
univariate polynomial $P^* \in R[T]$ and a non--zero element of the ring $\theta
\in R$, such that~:

        $$P^*=\theta {P\over gcd(P,P')}$$ holds. \end{lemma}

Suppose now that the characteristic of $k$, $K$ and $R$ is positive, say $p$. In
all our applications of the technical lemmas of this section $R$ will be a
polynomial ring over $k$ generated by some $k$--linear forms in the variables
$X_1,\ldots,X_n$, say $Y_1,\ldots,Y_r$ with $0 \le r \le n$ (see Section
\ref{alg}). As one easily verifies there is no essential change to the outcome
of our main algorithm in Section \ref{alg} if we replace the variables
$X_1,\ldots,X_n$ of our input equation system $f_1,\ldots,f_n$ by suitable
$p^k$--th powers of them with $k \in \N$ not too big (recall that we suppose the
ground field to be perfect and our arithmetic networks and straight--line
programs to include special gates for the extraction of $p$--th roots in $k$and
in ${\bar k}$). Thus we may suppose without loss of generality that we are able
to extract suitable $p$--th roots in
$R$ (and hence in $K$) of not too high degree ($p^k \le D$ will suffice).
Applying now  Lemma \ref{sqfree} iteratively (with at most $log_p D$ iterations)
we see that we obtain in sequential time $O(D^7)$ a separable representation
${\tilde P} \in R[T]$ of the polynomial $P$. Taking into account that in all our
applications the factor $\theta$ appearing in Lemma \ref{sqfree} will be monic
with respect to one of its variables and using standard tricks of parallelization
we get a division free and {\it well parallelizable} arithmetic network of size
$D^{O(1)}$ which from the coefficients of $P$ computes the coefficients of a
separable representation ${\tilde P} \in R[T]$ of $P$ (see \cite{giu-he92}, 2.1.
and \cite{gihesa} for more details about how to modify algorithms in order to
achieve the assumptions which allows to compute ${\tilde P}$).

\subsection{Correct Test Sequences and ``Vermeidung von Divisionen"}
\label{cuestores} The algorithms in \cite{figismi}, \cite{giu-he92},
\cite{gihesa}, \cite{krick-pardo1}, \cite{krick-pardo:CRAS} rely heavily on the
use of ``correct test sequences" (\cite{cucketal}, \cite{heschnorr},
\cite{krick-pardo1}) as a tool for deciding identity of polynomials given by
straight--line programs. Unlike the probabilistic identity tests
\cite{demillolipton}, \cite{iba-mor}, \cite{schw}, \cite{zip} the choice of a
suitable correct test sequence does not depend on the specific polynomials whose
identity has to be checked but only on the number of variables and the size of
the input circuit. The outcome are nonuniform deterministic or random algorithms
of a stronger type as those proposed in  \cite{schw} and
\cite{zip}.

 \begin{definition}
 Let ${\cal F}$ be a set of polynomials of $k[X_1,\dots,X_n]$ such that $0$
belongs to ${\cal F}$. Let ${\cal Q}$ be a subset of $k^n$. ${\cal Q}$ is called
a correct test sequence (or questor set) for ${\cal F}$ if for any $P
\in {\cal F}$ the following implication holds~:
 $$ P(x)=0 \;for \; all \; x \in {\cal Q} \Longrightarrow  P=0\ . $$

 \end{definition}

The cardinality $\# {\cal Q}$ is called the {\it length} of the test sequence
${\cal Q}$. 

The existence of short correct test sequences is warranted by the following
fact~:

\begin{lemma}
 \label{quest1}   Let ${\cal F}$ be the class of all polynomials of
$k[X_1,\ldots,X_n]$ which can be evaluated by a straight--line program in
${\bar k}(X_1,\ldots,X_n)$ of (nonscalar) size $L$ and depth $\ell$.  Let 
$\omega ~\!:=(2^{\ell +1}-2)\, (2^{\ell}
 +1)^2 \qquad {\rm and} \qquad
\sigma~\!:= 6 \,(\ell L)^2$. Then for any collection $\Omega$ of $\omega$
elements of $k$ the set $\Omega^n \subset k^n$ contains at least
$\omega^{n\sigma}\,(1-\omega^{-{ \sigma \over 6}})$ correct test sequences of
length $\sigma$ for the class ${\cal F}$.
 \end{lemma}

If the characteristic of $k$ is zero we can always make the standard choice
$\Omega\!:=\{1,\ldots,\omega \} \subset \Z$.

Lemma \ref{quest1} implies that for $L > 0$ and $\ell > 0$ there always exist
correct test sequences for ${\cal F}$ having length $\sigma$ which is polynomial
in $L$, and that any random choice of $\sigma$ elements in
$\Omega^n$ leads to a correct test sequence for ${\cal F}$ with an error
probability of
$\omega^{-{ \sigma \over 6}} << {1 \over 2}$. The proof of Lemma
\ref{quest1} can be found in  \cite{krick-pardo1} (see also \cite{heschnorr}).

We shall make frequent use of a method called ``Vermeidung von Divisionen" which
is due to \cite{stras:verm} and which allows to transform any straight--line
program $\Gamma$ (with parameters in $k$ or ${\bar k}$) which computes a
polynomial and contains essential divisions into an equivalent division free
one, say
$\Gamma '$. If the size and the depth of $\Gamma$ are $L$ and $\ell$
respectively, then the size and depth of $\Gamma'$ are $D^2L$ and $O(\ell)$ (see
\cite{krick-pardo1} for the proof of this precise statement).

\subsection{Construction of a Primitive Element}
\label{primitivo}

Primitive element constructions (or ``shape lemmas") are crucial for solving
zerodimensional polynomial systems symbolically (see \cite{royetal}, 
\cite{can}, \cite{chis}, \cite{chigri}, \cite{giu-he92}, \cite{lalo-lupe},
\cite{he-morg},
\cite{krick-pardo1}, \cite{lala}, \cite{moeller}). They furnish a handy
description of the multiplication tensor of the finite dimensional
$k-algebras$ which arise as (typically reduced) coordinate rings of such
equation systems.

Let $R$ be an integral $k$--algebra with fraction field $K$. Let
$I$ be a zero-dimensional ideal of $K[X_1,\ldots,X_n]$. We denote the radical of
$I$ by $\sqrt I$.

 For any  maximal ideal ${\cal Q}$ associated to $I$, let $K({\cal Q})$ be the
field 
$$K({\cal Q})\!:= K[X_1,\ldots,X_n]/{\cal Q}.$$
 A linear form 
$u=t_1X_1+\cdots+t_nX_n \in R[X_1,\ldots,X_n]$ is said to be a {\sl primitive
element} for $I$ (or for the finite dimensional $K-algebra$
$K[X_1,\ldots,X_n]/{\sqrt I}$) if it satisfies the following two conditions~: 

\begin{itemize}

\item for any maximal ideal ${\cal Q}$ associated to $I$, the minimal
      equation (minimal polynomial) of $u$ modulo ${\cal Q}$ has degree 
      $[K({\cal Q}): K]_{sep}$ (here $[K({\cal Q}): K]_{sep}$ denotes the
separability degree of the finite field extension $K \longrightarrow K({\cal
Q})$).

 \item For any two different maximal ideals ${\cal Q}$ and ${\cal
       Q}'$ associated to $I$, the minimal polynomials of $u$ modulo 
       ${\cal Q}$ and ${\cal Q}'$ are distinct.  

\end{itemize}

If the algebra $K[X_1,\ldots,X_n]/{\sqrt I}$ is unramified over $K$ then a linear
form $u \in R[X_1,\ldots,X_n]$ represents a primitive element for $I$ if and only
if the set
$\{1,u,u^2,u^3,\ldots\}$ generates modulo
$\sqrt I$ the $K-$vector space $K[X_1,\ldots,X_n]/{\sqrt I}$. We shall also apply
the  notion of ``primitive element" to zerodimensional localizations (by a given
polynomial) of generally nonzerodimensional ideals. When doing so (e.g. in Lemma
\ref{sinfortorico}) we have to think $I$ as the intersection of all primary
components of the given ideal which do not contain the polynomial.

For the construction of a primitive element from a given set of generators of
$I$ we follow the lines of \cite{giu-he92} or
\cite{krick-pardo1}.In order to explain our method we introduce new variables
$T_1,\ldots,T_n$ and for $1 \le j \le n$ we consider the following domains~: 

$$R_j~\!:=R[T_1,\ldots,T_{j-1},T_{j+1},\ldots,T_n]$$
$$K_j~\!:=K(T_1,\ldots,T_{j-1},T_{j+1},\ldots,T_n),$$ and the following linear
form in $X_1,\ldots,X_{j-1},X_{j+1},\ldots,X_n$~:

$$Z_j\!:=T_1X_1+\cdots+T_{j-1}X_{j-1}+T_{j+1}X_{j+1}+\cdots+T_nX_n.$$ 

Furthermore we consider the following linear form $U$ in $X_1,\ldots,X_n$~:
          $$U\!:=T_1X_1+\cdots+T_nX_n.$$  We observe that for any $1 \le j \le
n$ the identity $U=Z_j+T_jX_j$ holds. Let
$T$ be an additional variable and $D$ a natural number. The way how we shall
generate a primitive element is contained in the following statement~:

\begin{lemma}
\label{primitive} Let be given the following inputs~:

\begin{itemize}

\item monic polynomials $f_1,\ldots,f_n \in R[T]$ of degree at most $D$
      such that for each $1 \le j \le n$ the specialization $f_j(X_j)$ belongs to
      the ideal $I$.
\item for each $1 \le j \le n$ a polynomial $g_j\in R_j[T]$ which is monic in T
      and which has total degree at most $D$ such that the specialization
      $g_j(Z_j)$ belongs to the ideal $K_j\bigotimes_K I$.
\end{itemize}

We suppose further that $f_1,\ldots,f_n$ and $g_1,\ldots,g_n$ are represented by
their coefficients with respect to the variable $T$. Moreover we assume that the
coefficients of $g_1,\ldots,g_n$ themselves are given by a division free
straight--line program in $R[T_1,\ldots,T_n]$ of length $L$ and depth
$\ell$.

Then there exists a division free arithmetic network in $R$ of size
$(nDL)^{O(1)}$ and depth $O(log_2 (nD)+ \ell)$ which computes the coefficient
representation of the following items~:

\begin{enumerate}
\item a primitive element $u= \lambda_1X_1+\cdots+\lambda_n X_n \in
      k[X_1,\ldots,X_n]$ for the ideal $I$,
\item a non--zero element $\rho \in R$ and for each $1 \le j \le n$ a
       polynomial $v_j \in R[T]$ such that $\rho X_j-v_j(u)\in \sqrt {I}$ holds,
\item a monic polynomial $q\in R[T]$ such that $q(u)\in \sqrt {I}$
      holds.
\end{enumerate}

Moreover we have $deg v_i < degq$ for $1 \le i \le n$. From the observations at
the end of Section \ref{subsecuno} we deduce that we may suppose  without loss
of generality that
$q$ is separable with respect to the variable $T$ and hence squarefree. This
implies that the algebra $K[X_1,\ldots,X_n]/{\sqrt I}$ is unramified over $K$ and
that ${\sqrt I}=(q(u),\rho X_1-v_1(u),\ldots,\rho X_n-v_n(u))$ holds in
$K[X_1,\ldots,X_n]$.

\end{lemma}

Lemma \ref{primitive} is contained in \cite{krick-pardo1}, Proposition 27 and
its proof (see also \cite{giu-he92}, Sections 3.4.6 and 3.4.7).

\section{The Algorithm}
\label{alg} In this Section we describe the algorithmic procedure underlying
Theorem
\ref{afin}, \ref{tor} and \ref{parametrizacion}. 

Let us fix the following notions and notations~: let $X_0$ be a new variable.
For any non--zero polynomial $p\in k[X_1,\ldots,X_n]$  of degree $D$ we define
its homogenization $^h\!p$ (with  respect to the variable $X_0$) as 
$$^h\!p\!:= X_0^D p(X_1/X_0,\ldots,X_n/X_0).$$ For $p=0$ we put $^h p\!:=0$. Of
course $^h p$ is a homogeneous polynomial in $k[X_0,\ldots,X_n]$ of degree
$D$.  For any ideal $ I$ of $k[X_1,\ldots,X_n]$,  we define  its homogenization
$^h I$ as  the (homogeneous) ideal in $k[X_0,\ldots,X_n]$ generated by the set
of polynomials
$\{^h\!f~: f\in I\}$.

We maintain assumptions and notations of Section
\ref{introduccion}. In particular we suppose to be given input polynomials
$f_1,\ldots,f_n \in k[X_1,\ldots,X_n]$ of degree at most $d$ which are encoded
by a division free straight--line program of size $L$ and depth $\ell$. We
suppose that
$f_1,\ldots,f_n$ form an affine or toric complete intersection, according to the
problem we are considering.  Let
$H$ be a non--zero linear form of $k[X_1,\ldots,X_n]$ given by its coefficients.
The fundamental problem we want to solve is the following~:  find a non--zero
polynomial
$p$ (in the affine case) or a non--zero polynomial $p^{*}$ (in the toric case)
belonging to $k[T]$ such that $p(H)$ vanishes on
$V=V(f_1,\ldots,f_n)$ (in the affine case) or $p^{*}(H)$ vanishes on 
$V^{*}=V(f_1,\ldots,f_n) \setminus V(\prod\limits_{i=1}^n X_i)$ (in the toric
case). This is the content of Theorem \ref{afin} and 
\ref{tor}.

>From \cite{krick-pardo1} Lemma 13 we deduce that the homogeneous polynomials  
$^h\!f_1$,$\ldots$, $^h\!f_n$ can be evaluated  by a division free
straight--line program in $k[X_0,\ldots,X_n]$ of (nonscalar) size
 and depth $d(d+1)^2L$ and $log_2 d+2\ell$ respectively.

For $1 \le i \le n$ we denote by $I_i\!:=(f_1,\ldots,f_i)$ the ideal generated by
$f_1,\ldots,f_i$ in
$k[X_1,\ldots,X_n]$. Let us now introduce the specific notions and notations
which we need in the {\it affine case}. Fix $1 \le i \le n$ and let $
J_i~\!:=\sqrt {^h I_i}$ be the radical of the homogenization of
$I_i=(f_1,\ldots,f_i)$. Thus
$J_i$ is a homogeneous and unmixed radical ideal of codimension $i$ and the
corresponding projective variety does not contain any irreducible component at
infinity. Furthermore the homogeneous polynomial
$^h\!f_{i+1}$ is a non--zero divisor modulo $J_i$ for
$i<n$. Let ${\cal I}_i\!:=(^h\!f_i,J_{i-1})$ be the homogeneous ideal generated
by $^h\!f_i$ and
$J_{i-1}$. We observe that
$J_i$ is the intersection of all codimension
$i$ prime homogeneous ideals which contain ${\cal I}_i$ but not the form $X_0$
(of course, these are associated primes of ${\cal I}_i$). 

Let us finally consider the specific notions and notations of the {\it toric
case}. Let $1 \le i \le n$.  We denote by
$J_i^{*}$ the intersection of  all homogeneous prime ideals of codimension $i$
which contain
$^h I_i$ but not the form $\prod\limits_{i=0}^n X_i$ (these prime ideals are
again associated to $^h I_i$). By hypothesis (the family $f_1,\ldots,f_n$ forms
a toric complete intersection) the variety $V(I_i)$ is not contained in
$V(\prod\limits_{i=1}^n X_i)$ and therefore $ J_i^* \ne (X_0,\ldots,X_n)$.
Furthermore the homogeneous polynomial
$^h\!f_{i+1}$ is not a zero divisor modulo $J_i^*$ for $i < n$. Finally, as in
the affine case, $J_i^{*}$ is the intersection of all codimension $i$
homogeneous (associated) prime ideals which contain
$(^h\!f_i,J_{i-1}^*)$ but not the form
$\prod\limits_{i=0}^n X_i$.

With these notions and notations fixed we are now able to describe the principal
items which are produced as intermediate results by the algorithm underlying
Theorems
\ref{afin}, \ref{tor} and \ref{parametrizacion} and to indicate their main
properties.

For the affine as well as for the toric case the algorithm proceeds in an
analogous manner in $n$ recursive steps. We give now a simultaneous account  of 
the $i$--step of the  algorithm both for the affine and the toric case (here $1
\le i \le n$).  The algorithm produces first a linear change of the variables 
$X_0,\ldots,X_n$ into new variables $Y_0,\ldots,Y_n$ such that the following
canonical ring homomorphism becomes a generically unramified integral 
extension~:

\begin{itemize} 

\item {\sl in the affine case~:}
      $$k[Y_0,\ldots,Y_{n-i}] \longrightarrow k[Y_0,\ldots,Y_n]/J_i$$

\item {\sl in the toric case~:}
       $$k[Y_0,\ldots,Y_{n-i}] \longrightarrow k[Y_0,\ldots,Y_n]/J_i^{*}$$
 
\end{itemize}

Next the algorithm generates a nonzero linear form 
$u_i \in k[Y_{n-i+1},\ldots,Y_n]$ and a homogeneous polynomial of
$k[Y_0,\ldots,Y_{n-i},T]$ which is monic and separable with respect to the
variable $T$ and which we denote by 
$q_i $ in the affine and by  $q_i^{*}$ in the toric case. The linear form $u_i$
is a primitive element for the ideal
$J_i$ in the affine and for $J_i^{*}$ in the toric case having minimal
polynomial $q_i$ or
$q_i^{*}$ respectively. Finally the algorithm produces in the affine case
homogeneous polynomials
$\rho_i \in k[Y_0,\ldots,Y_{n-i}]$ and
$v_{n-i+1}^{(i)},\ldots,v_n^{(i)} \in k[Y_0,\ldots,Y_{n-i},T]$ with $\rho_i \ne
0$ such that in the localized ring
$k[Y_0,\ldots,Y_{n-i}]_{\rho_i}[Y_{n-i+1}\ldots,Y_n]$ the ideal generated by 

   $$q_i(u_i),\rho_i Y_{n-i+1} - v_{n-i+1}^{(i)}(u_i),\ldots, \rho_i Y_n-
     v_n^{(i)}(u_i)$$ is identical with the ideal $(J_i)_{\rho}$. The polynomial
$\rho_i$ will be the discriminant of the (separable) polynomial $q_i$ with
respect to the variable $T$. 

A similar statement is valid for the toric case if we replace
$q_i$ by $q_i^{*}$ and $J_i$ by
$J_i^{*}$. In order to avoid accumulation of diacritic symbols in the proof we
shall not use extra notation for the polynomials which play in the toric case
the r\^ole of $\rho_i$ and $v_{n-i+1}^{(i)},\ldots,v_n^{(i)}$. The same notation
will be applied in the affine case as well as in the toric case, no ambiguity
will arise from that. Thus for example $\rho_i$ will be interpreted as the
discriminant of $q_i$ or as the discriminant of $q_i^{*}$, following the context.

Let us also remark that the degree of $v_{n-i+1}^{(i)},\ldots,v_n^{(i)}$ with
respect to the variable $T$ will be strictly less than the degree of $q_i$ or
$q_i^{*}$ respectively (this latter condition makes
$v_{n-i+1}^{(i)},\ldots,v_n^{(i)}$ unique in
$k(Y_0,\ldots,Y_{n-i})[T]$).

Thus the $k(Y_0,\ldots,Y_{n-i})-$algebra
$k(Y_0,\ldots,Y_{n-i})[T]/(q_i)$ (or equivalently

\noindent $k(Y_0,\ldots,Y_{n-i})[u_i]/{(q_i(u_i))}$) is isomorphic to 
$k(Y_0,\ldots,Y_{n-i})[Y_{n-i+1},\ldots,Y_n]/(J_i)$ in the affine case (and the
isomorphism can be effectively computed). Both algebras are unramified. An
analogous statement is true in the toric case. 

The linear form $u_i$ will be given by its coefficients with respect  to the
variables $Y_{n-i+1},\ldots,Y_n$, whereas $q_i$,
$q_i^{*}$, $v_{n-i+1}^{(i)},\ldots,v_n^{(i)}$ will be given by their
coefficients (which belong to $k[Y_0,\ldots,Y_{n-i}]$) with respect to the
variable $T$. These coefficients and the polynomial $\rho_i$ will be represented
by a division free straight--line program in ${\bar k}[Y_0,\ldots,Y_{n-i}]$).
Furthermore our procedure will produce this straight--line program.

Finally let us observe that these items which appear as intermediate results of
our procedure are canonical and intrinsic objects with a precise geometric
meaning (on just this meaning is based the data compression contained in Lemma
\ref{sinfor} and
\ref{sinfortorico} which is fundamental for our approach). Let us consider only
the affine case~: to the linear change of the variables
$X_0,\ldots,X_n$ into $Y_0,\ldots,Y_n$ corresponds a finite surjective morphism
of affine varieties $\pi_i: V(f_1,\ldots,f_i)
\longrightarrow {\bf A}^{n-i}$. The linear form $u_i$  complements $\pi_i$ to a
finite morphism $(\pi_i,u_i)$ which maps $V(f_1,\ldots,f_i)$ onto a hypersurface
of ${\bf A}^{n-i+1}$. The minimal equation of this hypersurface is just the
polynomial we obtain specializing in $q_i \in k[X_0,\ldots,X_n]$ the variable
$X_0$ to $1$. With the same type of specialization we see that $\rho_i$
describes the discriminant variety of the generic fiber of $\pi_i$ and that the
polynomials
$v_{n-i}^{(i)},\ldots, v_n^{(i)}$ represent a ``univariate" rational
parameterization of the $\pi_i$--fibers of the points of ${\bf A}^{n-i}$ which 
lie outside of this discriminant variety.

\subsection{The Recursion}
\label{induc}

Keep $ 1 \le i < n$ and notions and notations of the last subsection fixed. We
are going to describe the $i+1$--th step of our recursive main algorithm. Let
$\delta_i\!:=degV(f_1,\ldots,f_i)$ and
$\delta_i^{*}\!:=deg^{*}V(f_1,\ldots,f_i)$ and recall that
$deg f_j \le d$ holds for $1 \le j \le n$. With these notations we have
$deg_T q_i=deg q_i \le \delta_i$,
$deg_T q_i^{*}=deg q_i^{*} \le \delta_i^{*}$.  {\it In order to simplify
notations we shall assume from now on that the polynomials $q_i$ and
$q_i^{*}$ have degree exactly $\delta_i$ and $\delta_i^{*}$ respectively. This
assumption does not restrict the generality of our arguments (compare Section
\ref{compresion}).} Let
$\Gamma_i$ be a division free straight--line program in
${\bar k}[Y_0,\ldots,Y_{n-i}]$ which  represents the coefficients of $u_i$ and
which evaluates $\rho_i$ and the coefficients with respect to $T$ of the
polynomials $v_{n-i+1}^{(i)},\ldots,v_n^{(i)}$ and of $q_i$ in the affine, and
of $q_i^{*}$ in the toric case. Denote by
$\Lambda_i$ and $\lambda_i$ the size and depth of $\Gamma_i$ respectively.

The core of our procedure is the following technical result~:

\begin{proposition} 
\label{inductive} There exists an arithmetic network with parameters in ${\bar
k}$ of 
\begin{itemize}
\item size $i(d\delta_i L\Lambda_i)^{O(1)}$ 
\item depth
$O(log_2(d\delta_i)+\ell)+\lambda_i$ 
\end{itemize} which from
$\Gamma_i$ as input produces a linear change of variables
$Y_0,\ldots,Y_n$ into new variables, say
$Y_0',\ldots,Y_n'$, such that
$Y_0',\ldots,Y_n'$ are in Noether position with respect to the ideal $J_{i+1}$
in the affine case and to $J^{*}_{i+1}$ in the toric case, the variables
$Y'_0,\ldots,Y'_{n-i-1}$ being free. 

Furthermore the network produces the coefficients of a linear form $u_{i+1}$ in
the ring $k[Y'_{n-i},\ldots,Y'_n]$ and a division free straight--line program in
${\bar k}[Y'_0,\ldots,Y'_{n-i-1}]$ which represents the polynomial
$\rho_{i+1}$ and the coefficients with respect to $T$ of
 $v_{n-i}^{({i+1})},\ldots,v_n^{({i+1})} \in k[Y'_0,\ldots,Y'_{n-i-1},T]$ and of
$q_{i+1}$ or $q_{i+1}^{*} \in k[Y'_0,\ldots,Y'_{n-i-1},T]$,  following the
(affine or toric) case. 

The size and depth of this straight--line program are
$i(L+\Lambda_i)(d\delta_i)^{O(1)}$ and
$O(log_2(d\delta_i)+\ell)+\lambda_i$ respectively. The parameters of the
arithmetic network and of the straight--line program it produces are contained
in the field extension of $k$ generated by the parameters of $\Gamma_i$.
\end{proposition}

We observe without proof that the algorithm underlying this procedure can be
organized in such a way that it uses only linear algebra subroutines dealing
with square matrices of size at most
$2d
\delta_i$ (or 
$2d \delta_i^{*}$).

\bigskip
 
We divide the $i+1$--th recursive step of our algorithm in three parts~: 

\begin{itemize}

\item {\sl recursive Noether normalization.} 
\item {\sl recursive generation of a primitive element.} 
\item {\sl cleaning extraneous irreducible components.}

\end{itemize}

The first two parts are common for both the affine and the toric case. Only the
third  part distinguishes in some few technical points between the two cases.

Before giving the details of these three parts, let us recall the following
procedure to compute resultants which is implicitly contained in 
\cite{can-gal-he} and \cite{saso}~:

\medskip

\begin{lemma}
\label{resultante} Let $J$ be an unmixed homogeneous radical ideal of 
$k[X_0,\ldots,X_n]$ and let us suppose that the canonical homomorphism
$$\phi : k[X_0,\ldots,X_{n-i}] \longrightarrow k[X_0,\ldots,X_n]/J.$$ is an
integral graded ring extension. Then, for an homogeneous polynomial $g$ in the
ring $k[X_0,\ldots,X_n]$ its minimal polynomial $m_g \in
k[X_0,\ldots,X_{n-i},T]$ has the form
$$T^D+a_{D-1}T^{D-1}+\cdots+a_0,$$ where $ D\leq deg J $ and $a_j \in
k[X_0,\ldots,X_{n-i}]$ is zero or a homogeneous polynomial of degree
$(D-j)\cdot deg(g)$.

The same conclusion remains true if we replace the minimal polynomial $m_g$ by
the characteristic polynomial $\chi_g$ of $g$. Moreover the minimal polynomial
$m_g$ is squarefree if $\phi$ is generically unramified and we have $D\leq deg
V(J)$ in this case. 
\end{lemma}

{\it Proof.} Observe first that our assumptions on $J$ imply that
$m_g$ and $\chi_g$ belong to $k[X_0,\ldots,X_{n-i},T]$. 

Let $J= \bigcap\limits_{j=1}^N {\cal Q}_j$, where the
${\cal Q}_j$'s are the homogeneous prime ideals of codimension $i$ associated to
$J$. From our assumptions we obtain for $1 \le j \le N$ an integral ring
extension~:
$$ \phi_j:k[X_0,\ldots,X_{n-i}]\hookrightarrow k[X_0,\ldots,X_n]/{\cal Q}_j$$
which is generically unramified if the same holds for $\phi$. We consider the
multiplication by the polynomial $g$ modulo the ideal ${\cal Q}_j$ as a
$k[X_0,\ldots,X_{n-i}]$--linear map. In this sense we shall speak in future
about ``the homothety defined by multiplication by $g$ modulo ${\cal Q}_j$".

\noi The minimal polynomial $m_g^{(j)}$ of this homothety verifies the
conclusions of the lemma with respect to degrees (see
\cite{can-gal-he}, Remark $9$ and its proof). Moreover $\phi_j$ generically
unramified means that the corresponding finite field extension is separable.
Therefore the minimal polynomial $m_g$ is squarefree if
$\phi$ is generically unramified. Since the minimal polynomial of the homothety
defined by the multiplication by
$g$ modulo
$J$ is a product of some of the given
$m_g^{(1)},\ldots,m_g^{(N)}$ (without repetitions if $\phi$ is generically
unramified), the first assertion of the lemma follows. The third one is a
consequence of the fact that the polynomials
$m_g^{(1)},\ldots,m_g^{(N)}$ are irreducible and that $J$ is radical. The second
assertion follows from the observation that the irreducible factors of
$\chi_g$ are the polynomials $m_g^{(1)},\ldots,m_g^{(N)}$.
\cqfd

Note that if $g$ is not a zero divisor modulo $J$, the term $a_0$ is a non--zero
homogeneous polynomial of $k[X_0,\ldots,X_{n-i}]$ of degree $D\cdot deg(g)$.

\subsubsection{Recursive Noether Normalization.} 
\label{uno} 

Without loss of generality we describe just the affine case (the toric case
follows simply by replacing
$J_i$ by 
$J_i^{*}$ and so on).

First we observe that $^h\!f_{i+1}$ is not a zero divisor modulo
$J_i$. According to the notations of Lemma \ref{resultante} above let
$\displaystyle{m_{^h\!f_{i+1}}}$ be the minimal polynomial of the homothety
given  modulo $J_i$ by $^h\!f_{i+1}$  ~:
$$m_{^h\!f_{i+1}}=T^D +a_{D-1} T^{D-1} + \cdots + a_0.$$ The coefficient $a_0
\in k[Y_0,\ldots,Y_{n-i}]$ is a non-zero homogeneous polynomial of degree at
most $d \delta_i$ and
$\displaystyle{m_{^h\!f_{i+1}}}$ is squarefree by Lemma
\ref{resultante}.

In order to obtain a Noether normalization with respect to 
${\cal I}_{i+1}=(^h f_{i+1}, J_i)$ (and hence with respect to $J_{i+1}$) we have
just to find a $k$--linear change of the variables
$ Y_0,\ldots,Y_{n-i}$ into new ones, say  $Y_0',
\ldots,Y_{n-i}'$, such that the (uniquely defined) polynomial
$A_0(Y_0',\ldots,Y_{n-i}') \in k[Y'_0,\ldots,Y'_{n-i}]=k[Y_0,\ldots,Y_{n-i}]$
which satisfies the equality
$$A_0(Y_0',\ldots,Y_{n-i}') = a_0(Y_0,\ldots,Y_{n-i})$$ is monic in the variable
$Y'_{n-i}$.  We obtain this linear change of variables as follows~:
 
First observe that modulo $J_i$ the forms
$1,u_i,\ldots,u_i^{\delta_i -1}$ represent a basis of the

\noindent $k(Y_0,\ldots,Y_{n-i})$--vector space
$k(Y_0,\ldots,Y_{n-i})[Y_{n-i+1},\ldots,Y_n]/(J_i)$. By assumption the
coefficients with respect to the variable
$T$ of the polynomials 
$q_i$, $v_{n-i+1}^{(i)},\ldots,v_n^{(i)} \in k[Y_0,\ldots,Y_{n-i},T]$ as well as
the polynomial
$\rho_i \in k[Y_0,\ldots,Y_n]$ are given by the circuit
$\Gamma_i$. This implies that we can write down directly the (homogeneous)
entries of the companion matrix $M_i \in k[Y_0,\ldots,Y_{n-i}]^{\delta_i \times
\delta_i}$ of
$q_i$. The matrix $M_i$ describes the homothety defined by the multiplication by
$u_i$ in
$k(Y_0,\ldots,Y_{n-i})[Y_{n-i+1},\ldots,Y_n]/(J_i)$ with respect to the basis
given by the forms
$1,u_i,\ldots,u_i^{\delta_i-1}$. With respect to the same basis the homothety
defined by the multiplication by
$^h\!f_{i+1}$ has a matrix representation $M_{^h\!f_{i+1}}$ which satisfies the
equation~:
$$\rho_i^k M_{^h\!f_{i+1}}=
{^h\!f_{i+1}}(Y_0,\ldots,Y_{n-i},v_{n-i+1}^{(i)}(M_i),
\ldots,v_{n}^{(i)}(M_i) )$$ for some $k\leq d$ (in a slight abuse of notation we
write
$^h\!f_{i+1}(Y_0,\ldots,Y_n)$ for the polynomial obtained from the original
$^h\!f_{i+1}(X_0,\ldots,X_n)$ by means of the change of variables from 
$X_0,\ldots,X_n$ to $Y_0,\ldots,Y_n$ and multiplying in this new expression
submonomials in $Y_0,\ldots,Y_{n-i}$ by suitable powers of $\rho_i$).  Observe
that
$\rho_i^k M_{^h\!f_{i+1}}
\in k[Y_0,\ldots,Y_{n-i}]^{\delta_i \times \delta_i}$ holds and that the entries
are homogeneous polynomials. Therefore the entries of the matrix
$\rho_i^k M_{^h\!f_{i+1}}$ can be computed by an (ordinary) straight--line
program in
$k[Y_0,\ldots,Y_{n-i}]$ of length
$L\delta_i^{O(1)}$ and depth $O(log_2\delta_i)+ \ell$ from the coefficients of
$q_i$, $v_{n-i+1}^{(i)},\ldots,v_n^{(i)}$ and from
$\rho_i$. This is simply done by applying the evaluation
scheme for $^h\!f_{i+1}$ (which is given by the input)
to suitable entries of the matrices $v_{n-i+1}^{(i)}(M_i),
\ldots,v_{n}^{(i)}(M_i)$ which replace the variables $Y_{n-i},\ldots,Y_n$.

Next, we use Berkowitz's well parallelizable and division free algorithm
\cite{berk} in order to compute the (homogeneous) coefficients of the
characteristic polynomial
$\chi_{\rho_{i}^{k} {^h}\!f_{i+1}} \in k[Y_0,\ldots,Y_{n-i},T]$ of the homothety
given by $\rho_{i}^{k} {^h}\!f_{i+1}$ modulo $J_i$. The computation of the
coefficients of
$\chi_{\rho_{i}^{k} {^h}\!f_{i+1}}$ from the entries of
$\rho_i^k M_{^h\!f_{i+1}}$  requires an additional division free and well
parallelizable straight--line program in
$k[Y_0,\ldots,Y_{n-i}]$ of size $\delta_i^{O(1)}$.

Now taking into account our observations at the end of Section
\ref{subsecuno} we apply Lemma
\ref{sqfree} in order to obtain the (homogeneous) coefficients with respect to
$T$ of a separable representation
${\tilde \chi}_{\rho_{i}^{k} {^h}\!f_{i+1}} \in k[Y_0,\ldots,Y_{n-i},T]$ of
$\chi_{\rho_{i}^{k} {^h}\!f_{i+1}}$. This separable representation is without
loss of generality the minimal polynomial
$m_{\rho_{i}^{k} {^h}\!f_{i+1}}$ multiplied by some non--zero homogeneous
$\theta \in k[Y_0,\ldots,Y_{n-i}]$. Observe that the coefficients of
$m_{\rho_{i}^{k} {^h}\!f_{i+1}}$ are homogeneous too. The computation of
$\theta$ and the remaining coefficients of
${\tilde \chi}_{\rho_{i}^{k} {^h}\!f_{i+1}}=
\theta
\cdot m_{\rho_{i}^{k} {^h}\!f_{i+1}}$ requires one more division free and well
parallelizable straight--line program in $k[Y_0,\ldots,Y_{n-i}]$ of size
$\delta_i^{O(1)}$.

Putting everything together we obtain a division free straight--line program in
$${\bar k}[Y_0,\ldots,Y_{n-i}]$$
 of size $L\delta_i^{O(1)} +
\Lambda_i$ and depth $ O(log_2\delta_i)+ \ell+ \lambda_i $ which computes the
homogeneous constant coefficient, say
${\cal A} \in k[Y_0,\ldots,Y_{n-i}]$, of
$\theta \cdot m_{\rho_{i}^{k} {^h}\!f_{i+1}}$ and the polynomial $\theta$. We
remark that
${\cal A}$ has degree at most $(d\delta_i)^{O(1)}$. By Lemma \ref{resultante}
this coefficient has the form
${\cal A}=
\theta
\rho_i^m a_0$ for some
$m
\le d
\delta_i$  and is therefore nonzero. Since
$m_{^h\!f_{i+1}}$ annihilates
$^h\!f_{i+1}$ modulo $J_i$ we conclude that  $a_0 \in (^h f_{i+1},J_i)$ holds in
$k[Y_0,\ldots,Y_n]=k[X_0,\ldots,X_n]$. Since the polynomial
${\cal A}$ can be evaluated by a division free straight--line program in
${\bar k}[Y_0,\ldots,Y_{n-i}]$ of size $L(d \delta_i)^{O(1)} +
\Lambda_i$ and depth $O(log_2d\delta_i)+ \ell+\lambda_i$, there exists by Lemma
\ref{quest1} a correct test sequence of length $(L
\delta_i^{O(1)} +
\Lambda_i)^3$ in $k^{n-i+1}$ for this complexity class. Because the polynomial
${\cal A}$ is nonzero we find in sequential time $(L (d\delta_i)^{O(1)} +
\Lambda_i)^4$ and parallel time $O(log_2d\delta_i)+ \ell +\lambda_i$ a point
$\gamma=(\gamma_0,\ldots,\gamma_{n-i}) \in k^{n-i+1}$ in this correct test
sequence such that ${\cal A}(\gamma) \ne 0$ holds. Since ${\cal A}$ is
homogeneous we may construct in the obvious way from the coordinates
$\gamma_0,\ldots,\gamma_{n-i}$ of $\gamma$ a linear change of the variables
$Y_0,\ldots,Y_{n-i}$ into new ones, say $Y'_0,\ldots,Y'_{n-i}$, such that the
homogeneous polynomial ${\cal A}' \in
k[Y'_0,\ldots,Y'_{n-i}]=k[Y_0,\ldots,Y_{n-i}]$ given by the equation ${\cal
A}'(Y'_0,\ldots,Y'_{n-i})={\cal A}(Y_0,\ldots,Y_{n-i})$ becomes monic in
$Y'_{n-i}$. This implies that the homogeneous polynomial
$A_0(Y'_0,\ldots,Y'_{n-i})\!:=a_0(Y_0,\ldots,Y_{n-i})$ is monic too in
$Y'_{n-i}$.

Therefore, if we replace the variables $Y_0,\ldots,Y_n$ by the new variables 
$$Y'_0,\ldots,Y'_{n-i}, Y'_{n-i+1}\!:=Y_{n-i+1},\ldots,Y'_n\!:=Y_n$$
 we obtain a Noether normalization with respect to the ideal ${\cal
I}_{i+1}=(^hf_{i+1},J_i)$ and hence with respect to $J_{i+1}$. 

In the next part of our algorithm we shall make use of the eliminating form
$a_0$ instead of the ``accidental" polynomial ${\cal A}$ which comes from our
specific algorithm. We consider $a_0$,
${\cal A}$, $\theta$ and $\rho_i$ as forms in the variables
$Y'_0,\ldots,Y'_{n-i}$ which are related by the identity

           $$a_0= {{\cal A} \over {\theta \rho_i^k}}.$$

Let $\mu\!:= \theta \rho_i^k$ and observe that $\mu$ is monic in $Y'_{n-i}$. The
exponent $k \le d$ is known and the forms ${\cal A}$, $\theta$ and $\rho_i$ (and
hence $\mu=
\theta \rho_i^k$) are represented in our algorithm by a {\it division free}
straight--line program $\Sigma$ in
${\bar k}[Y'_0,\ldots,Y'_{n-i}]={\bar k}[Y_0,\ldots,Y_{n-i}]$ which has size
$L(d\delta_i)^{O(1)}+\Lambda_i$ and depth $O(log_2 d\delta_i))+\ell+\lambda_i$.
The degrees of $\mu$ and
${\cal A}$ are of order $(d\delta_i)^{O(1)}$.

Using $N\!:=max\{deg \mu,deg {\cal A} \}=(d\delta_i)^{O(1)}$ parallel organized
calls  to the procedure $\Sigma$ we interpolate ${\cal A}$ and $\mu$ with
respect to the variable $Y'_{n-i}$ in $N$ many distinct points of $k$ and
compute their coefficients with respect to the variable $Y'_{n-i-1}$ (these
coefficients are polynomials which belong to
$k[Y'_0,\ldots,Y'_{n-i-1}]$). This can be achieved by a division free
straight--line program in ${\bar k}[Y'_0,\ldots,Y'_{n-i}]$ of size and depth
$(L+\Lambda_i)(d\delta_i)^{O(1)}$ and $O(log_2d\delta_i)+ \ell +\lambda_i$
respectively. Since the square matrix with entries in
$k[Y'_0,\ldots,Y'_{n-i}]$ which correspond to the division (with remainder) of
${\cal A}$ by $\mu$ in the principal ideal domain
$k(Y'_0,\ldots,Y'_{n-i-1})[Y'_{n-i}]$ is unimodular and since the identity
$a_0= {{\cal A} \over \mu}$ holds, we are able to compute the coefficients of
the polynomial $a_0$ with respect to
$Y'_{n-i}$ (and hence $a_0$ itself) from the circuit
$\Sigma$  in sequential and parallel time 
$(d\delta_i)^{O(1)}$ and
$O(log_2d\delta_i)$ respectively. Therefore we may suppose that $a_0$ (and its
coefficients with respect to
$Y'_{n-i}$) is given by a {\it division free} straight--line program in 
${\bar k}[Y'_0,\ldots,Y'_{n-i}]={\bar k}[Y_0,\ldots,Y_{n-i}]$  of size and depth
$(L+\Lambda_i)(d\delta_i)^{O(1)}$ and $O(log_2d\delta_i)+\ell+ \lambda_i$
respectively.

\subsubsection{Recursive Generation of a Primitive Element.} 
\label{dos} 

The next recursive step is the computation of a primitive element
$u_{i+1}$ modulo $J_{i+1}$. We assume without loss of generality that the
variables $Y_0,\ldots,Y_n$ are already in Noether position with respect to the
ideal ${\cal I}_{i+1}=(^h\!f_{i+1},J_i)$ (see ``recursive Noether
normalization"). So, we have the following integral ring extension
$$R~\!:= k[Y_0,\ldots,Y_{n-i-1}] \hookrightarrow
k[Y_0,\ldots,Y_n]/(^h\!f_{i+1},J_i)=k[Y_0,\ldots,Y_n]/({\cal I}_{i+1})$$ Let us
denote the fraction field of $R$ by $K$ (i.e.
$K\!:=k(Y_0,\ldots,Y_{n-i-1})$). It is clear that
${\cal I}_{i+1}$ generates a zero dimensional ideal in
$K[Y_{n-i},\ldots,Y_n]$. We are going to apply the method described in Section
\ref{primitivo}. For this purpose we introduce new variables, say 
$T_{n-i},\ldots,T_n$, and with respect to them for $n-i \le j \le n$ the rings
and fields
$R_j\!:=R[T_{n-i},\ldots,T_{j-1},T_{j+1},\ldots,T_n]$ and
$K_j\!:=K(T_{n-i},\ldots,T_{j-1},T_{j+1},\ldots,T_n)$ introduced in
\ref{primitivo}. For $n-i \le j \le n$ we are going to construct polynomials
$h_j\in R[T]$ and $g_j\in R_j[T]$ with the following properties~:
 
\begin{itemize}

\item $h_j\in R[T]$ is monic in T of degree at most $d\delta_i$ in $T$
      and the polynomial $h_j(Y_j)$ belongs to the ideal
      ${\cal I}_{i+1}$
\item $g_j\in R_j[T]$ is  monic in $T$ and has total degree at most 
      $d\delta_i$. The polynomial $g_j(Z_j)$ belongs to the ideal
      $K_j\bigotimes_R {\cal I}_{i+1}$
 
\end{itemize}

Our aim is to represent the coefficients of $h_j$ and $g_j$ with respect to
$T$ by an evaluation scheme.

With the notations and results of ``recursive Noether normalization" there is
already given an equation $a_0 \in R[Y_{n-i}]=k[Y_0,\ldots,Y_{n-i}]$ for the
integral dependence of
$Y_{n-i}$ over $R$ modulo the ideal ${\cal I}_{i+1}$. This polynomial $a_0$ is
given by a straight--line program in
${\bar k}Ê\otimes_k R[Y_{n-i}]={\bar k}[Y_0,\ldots,Y_{n-i}]$ of length
$(L+\Lambda_i)(d\delta_i)^{O(1)}$ and depth $ O(log_2(d\delta_i)+\ell)
+\lambda_i$ and it has degree 
$d\delta_i$. Moreover it is homogeneous with respect to the natural grading of
$R[Y_{n-i}]$. 

Fix  $n-i \le j \le n$ and consider the variable $Y_j$ and the linear form 
$Z_j\!:=T_{n-i}Y_{n-i}+\cdots+T_{j-1}Y_{j-1}+T_{j+1}Y_{j+1}+\cdots+T_nY_n
\in R_j[Y_{n-i},\ldots,Y_n]$ in the variables
$Y_{n-i},\ldots,Y_n$ (see Section \ref{primitivo}). To the homotheties given by
the forms $Y_j$ and $Z_j$ in the $\delta_i$--dimensional algebras
$K[Y_{n-i},\ldots,Y_n]/(J_i)$ and $K_j[Y_{n-i},\ldots,Y_n]/(J_i)$ and to the
bases of these algebras determined by $1,u_i,\ldots,u_i^{\delta_i -1}$
correspond matrices
$M_{Y_j} \in K^{\delta_i \times \delta_i}$ and $M_{Z_j} \in K_j^{\delta_i \times
\delta_i}$. The characteristic polynomials of these matrices induce polynomials
in the variables
$Y_{n-i}, Y_j$ and $Y_{n-i},Z_j$ respectively which have degree at most
$\delta_i$ in these indeterminates. Since the variables
$Y_0,\ldots,Y_n$ are in Noether position with respect to
$J_i$ the coefficients of these polynomials  belong to
$R$ and
$R_j$ respectively. These coefficients are homogeneous elements of degree at most
$\delta_i$ of their respective graded rings. Let us denote the polynomials
introduced in this way as  
$H_j \in R[Y_{n-i},Y_j]$ and $G_j \in R_j[Y_{n-i},Z_j]$. They have the following
properties~:

\begin{itemize} 

\item $H_j$ is monic in $Y_j$ and belongs to the ideal $J_i$

\item $G_j$ is monic in $Z_j$ and belongs to the ideal 
$K_j\otimes_R J_i$. 

\end{itemize}

These polynomials can be computed from the results of the circuit
$\Gamma_i$ as in ``recursive Noether normalization" by means of the division
free and well parallelizable algorithm \cite{berk} in sequential time
$\delta_i^{O(1)}$. This means that their coefficients with respect to the
variables
$Y_j$ and $Z_j$ are given by division free straight--line programs in the rings
${\bar k} \otimes_kR[Y_{n-i}]={\bar k}[Y_0,\ldots,Y_{n-i}]$ and
${\bar k}
\otimes_kR_j[Y_{n-i}]={\bar
k}[T_{n-i},\ldots,T_{j-1},T_{j+1},\ldots,T_n,Y_{n-i}]$ and that these
straight--line programs have length
$\delta_i^{O(1)}+ \Lambda_i$ and depth
$O(log_2\delta_i)+\lambda_i$. (Without going into details we remark here that
one has to use the same type of subroutine to eliminate divisions by a certain
polynomial which is monic in $Y_{n-i}$ as we did at the end of ``recursive
Noether normalization" when computing $a_0$).

For $n-i \le j \le n$ we are now going to construct the announced polynomials
$h_j \in R[T]$ and
$g_j \in R_j[T]$. For this purpose we need the coefficients with respect to
$Y_{n-i}$ of $a_0$ and the coefficients with respect to 
$Y_{n-i}$, $Y_j$ and $Y_{n-i}$, $Z_j$ of $H_j$ and $G_j$. This requires
interpolation of these polynomials in
$d\delta_i+1$ points of $k$ which we have to substitute for the variable
$Y_{n-i}$ (recall that $H_j$ and
$G_j$ are already given by their coefficients with respect to $Y_j$ and $Z_j$
respectively). For this interpolation we need $(2i+1)(d\delta_i+1)$ parallel
organized calls to the whole procedure
 which originate a total sequential time cost of
$i(L+\Lambda_i)(d\delta_i)^{O(1)}$ and a parallel time cost of $O(log_2
d\delta_i+\ell)+\lambda_i$.
 
Thus the coefficients with respect to the variable
$Y_{n-i}$ of $a_0$ and for
$n-i \le j \le n$ the coefficients with respect to
$Y_{n-i}$, $Y_j$ and
$Y_{n-i}$, $Z_j$ of the polynomials $H_j$ and $G_j$ are represented by division
free straight--line programs in
${\bar k} \otimes_k R={\bar k}[Y_0,\ldots,Y_{n-i-1}]$ and ${\bar k} \otimes_k
R_j={\bar k}[T_{n-i},\ldots,T_{j-1},T_{j+1},\ldots,T_n]$ which have length
$i(L+\Lambda_i)(d\delta_i)^{O(1)}$ and depth
$O(log_2(d\delta_i)+ \ell) +\lambda_i$. Once these coefficients are given we are
able to perform effectively computations in the graded $k$--algebras which for
$n-i \le j \le n$ are defined as follows~:

\begin{itemize}

\item $B_j~\!:=R[Y_{n-i},Y_j]/(a_0,H_j)$

\item $C_j~\!:=R_j[Y_{n-i},Z_j]/(a_0,G_j)$.

\end{itemize}

Observe that $B_j$ is a free $R$--module with basis represented by the set of
monomials 
      $$ \{Y_{n-i}^\alpha Y_j^\beta~: 0\leq \alpha < deg \; (a_0),
         0\leq \beta < deg \; (H_j)\}$$  and that $C_j$ is a free $R_j$--module
with basis represented by the set of monomials 
      $$ \{Y_{n-i}^\alpha Z_j^\beta~: 0\leq \alpha < deg \; (a_0),
       0\leq \beta < deg \; (G_j)\}.$$

The ranks of the free modules $B_j$ and $C_j$ are therefore bounded by 
$d\delta_i^2$. Since the coefficients of $a_0$ with respect to
$Y_{n-i}$ and the coefficients of $H_j$ with respect to
$Y_{n-i}$, $Y_j$ are given, we know effectively (with no extra cost) the
multiplication tensor of
$B_j$. In the same sense the multiplication tensor of $C_j$ is available. From
this data we compute the characteristic polynomial $\chi_{\scriptscriptstyle
Y_j} \in R[T]$ of the homothety in $B_j$ given by the multiplication by $Y_j$
and the characteristic polynomial $\chi_{\scriptscriptstyle Z_j} \in R_j[T]$ of
the homothety in $C_j$ given by $Z_j$.  We put $h_j \!:=
\chi_{\scriptscriptstyle Y_j}$ and
$g_j\!:=\chi_{\scriptscriptstyle Z_j}$. One verifies immediately that
$h_j$ and $g_j$ have the required properties (for this purpose observe that
$a_0
\in {\cal I}_{i+1}$, $H_j \in J_i
\subseteq {\cal I}_{i+1}$ and $G_j \in K_j\otimes_R J_i
\subseteq  K_j\otimes_R {\cal I}_{i+1}$ holds).

Applying Lemma \ref{primitive} and the observations at the end of Section
\ref{subsecuno} we obtain a linear form
$u_{i+1}$ which is a primitive element for the ideal generated by ${\cal
I}_{i+1}$ in $K[Y_{n-i},\ldots,Y_n]$, a homogeneous polynomial
$q \in R[T]$ and homogeneous polynomials $\rho \in R$ and
$v_{n-i},\ldots,$\\$ v_n
\in R[T]$ such that the following holds~: $q$ is monic and separable in $T$
(hence squarefree), $q(u_{i+1})$ belongs to the ideal $\sqrt{{\cal I}_{i+1}}$,
$\rho$ is nonzero and for $n-i \le j \le n$ we have $\rho Y_j- v_j(u_{i+1}) \in
\sqrt{{\cal I}_{i+1}}$. Observe that we may suppose without loss of generality
that $\rho$ is monic in one of its variables, say $Y_{n-i-1}$.

These items (and their homogeneous coefficients in
$R=k[Y_0,\ldots,Y_{n-i-1}])$ are represented by a division free straight-- line
program in ${\bar k}[Y_0,\ldots,Y_{n-i-1}]$ whose length and depth are
$i(L+\Lambda_i)(d\delta_i)^{O(1)}$ and $O(log_2(id\delta_i)+
\ell)+\lambda_i$ respectively.

\subsubsection{Clearing Extraneous Irreducible Components.}  

In ``recursive Noether normalization" we obtained a fairly explicit description
of the projective variety defined by the ideal ${\cal I}_{i+1}$ and of the
localization of this ideal in
$K[Y_{n-i}\ldots,Y_n]$. However this projective variety may contain irreducible
components at infinity (extraneous components) or components which count in
${\cal I}_{i+1}$ (or its localization) with higher multiplicities than one.

The situation in the toric case is the same with ``extraneous components at
infinity" replaced by ``extraneous components contained in the union of
hyperplanes $V(\prod\limits_{i=1}^n X_i)$".  We show now how extraneous
components and multiplicities which appear in a natural way when cutting by the
hypersurfaces
$V(^h\!f_1),\ldots,V(^h\!f_n)$, can be cleared out during the process. Let us
concentrate upon the affine case.

First  recall that
 $J_{i+1}$  is the intersection of all those homogeneous prime ideals of
codimension $i+1$ which contain
$^h\!f_{i+1}$ and $J_i$ but do not contain
$X_0$. We observe also that the (homogeneous) primes of codimension $i+1$ of
$R[Y_{n-i},\ldots,Y_n]=k[Y_0,\ldots,Y_n]$ which contain the ideal  
 $(^h\!f_1,\ldots,^h\!f_{i+1})$ but not $X_0$ correspond to the primes  in the
localized ring $K[Y_{n-i},\ldots,Y_n]$ which contain the ideal generated  by
${\cal I}_{i+1}$ but not $X_0$. We are going now to determine the irreducible
components at finite distance of the projective variety defined by
$^h\!f_1,\ldots,^h\!f_{i+1}$ (in other words defined by ${\cal I}_{i+1}$) doing
computations in
$K[u_{i+1}]$. These components characterize geometrically the ideal
$J_{i+1}$. For
$1
\le j
\le i+1$ let us write
$F_j\!:= ^h\!f_j(Y_0,\ldots,Y_{n-i-1},{v_{n-i}(u_{i+1}) \over {\rho}},\ldots,
{v_n(u_{i+1}) \over {\rho}})
\in K[u_{i+1}]$. Let
$h\in R[u_{i+1}]$ be a greatest common divisor of the polynomials
$F_1(u_{i+1}),\ldots,F_{i+1}(u_{i+1})$ and $q(u_{i+1})$ with respect to the
principal domain
$K[u_{i+1}]$.

Since $q(u_{i+1})$ is monic and separable in $u_{i+1}$ we may suppose without
loss of generality that the polynomial $h$ is monic and separable too in
$u_{i+1}$. Moreover we may suppose that $h$ is homogeneous with respect to the
grading of $R[u_{i+1}]$. The coefficients of the polynomial 
$h$ with respect to the variable $u_{i+1}$ can be computed from the coefficients
of $F_1,\ldots,F_{i+1}$, $q \in K[u_{i+1}]$ by a well parallelizable and
division free straight--line program of size
$(d\delta_i)^{O(1)}$. This means that these coefficients can be computed by a
division free straight--line program in
${\bar k}Ê\otimes_k R ={\bar k}[Y_0,\ldots,Y_{n-i}]$ of size
$i(L+\Lambda_i)(d\delta_i)^{O(1)}$ and depth
$O(log_2(d\delta_i)+\ell)+\lambda_i$ (here using the assumption that $\rho$ is
monic in $Y_{n-i-1}$, we apply again the same trick as at the end of ``recursive
Noether normalization"  in order to eliminate division by
$\rho$).

Observe now that the maximal ideals in
$K[Y_{n-i},\ldots,Y_n]$ which contain the polynomial
$h(u_{i+1})$ and for $n-i \le j \le n$ the polynomials
$\rho Y_j -v_j(u_{i+1})$, are in one to one correspondence with the codimension
$i+1$ prime ideals in
$R[Y_{n-i},\ldots,Y_n]$ which contain
$(^h\!f_1,\ldots,^h\!f_{i+1})$. Moreover $h(u_{i+1})$ is squarefree. 

Note that $X_0 \in k[X_0,\ldots,X_n]= k[Y_0,\ldots,Y_n]$ is a linear form in
$Y_0,\ldots,Y_n$. Substituting in this linear form for the variables
$Y_{n-i},\ldots,Y_n$ the polynomials 
$${v_{n-i}(u_{i+1})
\over {\rho}},\ldots,{v_n(u_{i+1}) \over {\rho}} \in K[u_{i+1}]$$ and clearing
the denominator
$\rho$ we obtain a representation in
$R[u_{i+1}]$ of the residue class of $\rho X_0$ modulo the ideal generated by
${\cal I}_{i+1}$ in
$K[Y_{n-i},\ldots,Y_n]$. Let $G \in R[u_{i+1}]$ be this representation. Now we
compute the greatest common divisor of 
$G(u_{i+1})$ and $h(u_{i+1})$ in the principal ideal domain
$K[u_{i+1}]$. This greatest common divisor is represented by a squarefree
polynomial $h_1\in R[u_{i+1}]$ which divides $h \in R[u_{i+1}]$. Observe that
$h_1$ is monic in
$u_{i+1}$ and homogeneous with respect to the grading of
$R[u_{i+1}]$. The maximal ideals of  
$K[Y_{n-i},\ldots,Y_n]$ which contain $h_1(u_{i+1})$ and for
$n-i \le j \le n$ the polynomials $\rho Y_j -v_j(u_{i+1})$, are in one to one
correspondence with the codimension
$i+1$ prime ideals of $R[Y_{n-i},\ldots,Y_n]$ which contain
$(X_0, ^h\!f_{i+1},\ldots,^h\!f_n)$.

Now it is easy to see that the maximal ideals in
$K[Y_{n-i},\ldots,Y_n]$ which contain  for $n-i \le j \le n$ the polynomials
$\rho Y_j -v_j(u_{i+1})$ and the homogeneous and monic in $u_{i+1}$ polynomial
$$q_{i+1}(u_{i+1})\!:= { h(u_{i+1}) \over h_1(u_{i+1})}$$ correspond  exactly to
the codimension $i+1$ prime ideals which contain the radical ideal $J_{i+1}$
(these are also the associated primes of $J_{i+1}$). Observe that
$q_{i+1}\!:=q_{i+1}(T)$ is monic and separable with respect to the variable
$T$. Finally  the polynomials
$v_{n-i}^{(i+1)},\ldots,v_n^{(i+1)}$ are defined (and computed) as the
remainders of the division in
$R[T]$ of
$v_{n-i}(T)\ldots,v_n(T)$ by the monic in
$T$ polynomial
$q_{i+1}(T)$ (note that this leaves the ideal in question unchanged). In this
way we obtain an explicit description of an isomorphism 
$ K[Y_{n-i},\ldots,Y_{n}]/(J_{i+1}) \cong K[u_{i+1}]/(q_{i+1}(u_{i+1})) \cong
K[T]/(q_{i+1}(T))$ as wanted. Since $h_1$ is monic in $u_{i+1}$ we are able to
compute by the same trick as before $q_{i+1}$ by a {\it division free}
straight--line program. Because $q_{i+1}$ is again monic in $T$ the same
argument applies to
$v_{n-i}^{(i+1)},\ldots,v_n^{(i+1)}$. Finally we observe that $\rho$ and the
coefficients with respect to the variable $T$ of $q(T)$ and
$v_{n-i}^{(i+1)}(T),\ldots,v_n^{(i+1)}(T)$ can be computed by a division free
straight--line program $\Gamma$in ${\bar k}[Y_0,\ldots,Y_{n-i-1}]$ of size
$i(L+\Lambda_i)(d\delta_i)^{O(1)}$ and depth
$O(log_2(d\delta_i)+\ell)+\lambda_i$. Since $q_{i+1}$ is separable and monic
with respect to $T$ its discriminant (which we denote by $\rho_{i+1}$) is
nonzero and belongs to
$R=k[Y_0,\ldots,Y_{n-i-1}]$. Since it is possible to compute $\rho_{i+1}$ by a
well parallelizable straight--line program in $R$ of size
$(d\delta_i)^{O(1)}$ we may suppose without loss of generality that the circuit
$\Gamma$ represents also the polynomial $\rho_{i+1}$. 

We conclude this subsection with the remark that the toric case can be treated
in exactly the same way replacing at the end of the construction for the affine
case the linear form
$X_0$ by the homogeneous polynomial $\prod\limits_{i=0}^n X_i$.

\subsection{Proof of Theorems 1, 2 and 3}

In this subsection we deduce from Lemma \ref{sinfor} and \ref{sinfortorico} and
from Proposition \ref{inductive} first Theorem \ref{parametrizacion} and then
Theorem \ref{afin} and \ref{tor}.

\medskip

{\it Proof of Theorem \ref{parametrizacion}}. Let assumptions and notations be
the same as in the statement of Theorem
\ref{parametrizacion} and as at the beginning of this section. The algorithm
underlying the proof proceeds in $n$ recursive steps. 

Let us now show Theorem
\ref{parametrizacion} $(i)$ which deals with the affine case. In the first step
of our algorithm we apply Proposition \ref{inductive} just to the input
straight--line program $\beta$ and the input polynomials $f_1,\ldots,f_n$ in the
following sense~: we put formally $i=0$ and assume $\Gamma_0$ to be the empty
circuit.

Suppose now $ 1 \le i < n$ and let us consider the $(i+1)$--th step of our
procedure. Let be given an arithmetic network ${\cal N}_i$ with parameters in
${\bar k}$ which produces the input for the $(i+1)$--th recursion step of our
procedure, namely a division free straight--line program
$\beta_i$ in
${\bar k}[Y_0,\ldots,Y_{n-i}]$ which evaluates $\rho_i$ and the coefficients
with respect to $T$ of the polynomials
$q_i$, $v_{n-i+1}^{(i)},\ldots,v_n^{(i)}$. Suppose  also that
$\beta_i$ contains the information about the coefficients of the linear form
$u_i$ which for example may be stored as parameters of
$\beta_i$. Assume inductively that the polynomial
$q_i$ is monic and separable with respect to the variable
$T$ and that $\rho_i$ is its discriminant. Let $L_i$ and
$\ell_i$ be the nonscalar size and depth of ${\cal N}_i$ respectively. Note that
the size and depth of the circuit $\beta_i$ do not exceed $L_i$ and
$\ell_i$ respectively. 

To the circuits $\beta$ and $\beta_i$, the polynomials $f_1,\ldots,f_i$,
$\rho_i$, $q_i$, $v_{n-i+1}^{(i)},\ldots,v_n^{(i)}$ and the $k$--linear form
$u_i$ we apply the compression algorithm underlying Lemma \ref{sinfor}. We
obtain a new division free straight--line program $\Gamma_i$ in ${\bar
k}[Y_0,\ldots,Y_{n-i}]$ which represents the coefficients of $u_i$, the
polynomial $\rho_i$ and the coefficients with respect to $T$ of $q_i$ and
$v_{n-i+1}^{(i)},\ldots,v_n^{(i)}$. From Lemma \ref{sinfor} we deduce that the
circuit $\Gamma_i$ has size and depth $O((i^5+L)\delta_i^{11})$ and $O((log_2 i
+ \ell) log_2 \delta _i)$.

We apply now to the circuit $\beta$ and $\Gamma_i$, the polynomials
$$f_1,\ldots,f_{i+1},\rho_i,q_i,v_{n-i+1}^{(i)},\ldots,v_n^{(i)}$$ and the
$k$--linear form $u_i$ the elimination procedure underlying Proposition
\ref{inductive}. We obtain a linear change of the variables $Y_0,\ldots,Y_n$ of
the $i$--th step into new variables $Y'_0,\ldots,Y'_n$ such that
$Y'_0,\ldots,Y'_n$ are in Noether position with respect to the ideal
$J_{i+1}$, the variables $Y'_0,\ldots,Y'_{n-i-1}$ being free. Moreover we obtain
a new $k$--linear form $u_{i+1}$ and new polynomials $\rho_{i+1} \in
k[Y'_0,\ldots,Y'_{n-i-1}]$ and $q_{i+1}$,
$v_{n-i}^{(i+1)},\ldots,v_n^{(i+1)}\in k[Y'_0,\ldots,Y'_{n-i-1},T]$. All these
items are represented by a division free straight--line program $\beta_{i+1}$ in
$k[Y'_0,\ldots,$ $\ldots,Y'_{n-i-1}]$ which is generated by an arithmetic network
${\cal N}_{i+1}$ with parameters in ${\bar k}$. From Proposition
\ref{inductive} we deduce that the size and depth of $\beta_{i+1}$ are
$O(i(L+(i^5+L)\delta_i^{11})(d\delta_i)^{O(1)})=i^6 L(d\delta_i)^{O(1)}$ and
$O((log_2(d\delta_i) + \ell) + (log_2i +\ell)log_2\delta_i)= O(log_2 (id) +\ell)
log_2 \delta_i)$.

Similarly the size $L_{i+1}$ and the depth $\ell_{i+1}$ of ${\cal N}_i$ are
$(iLd\delta_i)^{O(1)}+L_i$ and  $O(log_2 (id) +\ell) log_2 \delta_i)+\ell_i$.

Finally we obtain for $i\!:=n-1$ the output of the algorithm underlying Theorem
\ref{parametrizacion}$(i)$. This output consists in the coefficients of a 
nonzero linear form $u \in k[X_1,\ldots,X_n]$ and of one--variate polynomials
$q$, $v_1,\ldots,v_n \in k[u]$ such that the following holds~:
\begin{itemize}
\item $deg q = \delta_n$
\item $max \{deg v_1,\ldots,deg v_n \} < \delta_n$
\item $(f_1,\ldots,f_n) = (q,X_1-v_1,\ldots,,X_n-v_n).$
\end{itemize}

The output is represented by the arithmetic network ${\cal N}_n$ whose
parameters belong to ${\bar k}$ and whose size and depth are
$(nLd\delta)^{O(1)}$ and $O(n(log_2(nd)+\ell)log_2\delta)$ respectively.

The proof of Theorem \ref{parametrizacion} $(ii)$ which deals with the toric
case is completely analogous. One has just to replace in it the application of
Lemma \ref{sinfortorico} by the application of Lemma \ref{sinfor} in the same
sense. In order to avoid repetitive arguments we omit this proof.
\cqfd

\medskip

>From Theorem \ref{parametrizacion} $(i)$ and $(ii)$ we deduce now easily Theorem
\ref{afin} and \ref{tor}.

\medskip

{\it Proof of Theorem \ref{afin} and \ref{tor}.} Let us first concentrate upon
the affine case, namely Theorem \ref{afin}.

Let assumptions and notations be the same as in the statement of Theorem
\ref{afin}. Applying Theorem \ref{parametrizacion} $(i)$ to this situation we
obtain an arithmetic network ${\cal N}$ with parameters in ${\bar k}$ which
computes the coefficients of a $k$--linear form $u$ and of polynomials $q$,
$v_1,\ldots,v_n
\in k[T]$ satisfying the following conditions~:

\medskip

(*) $deg q= \delta_n$ and $max \{deg v_i \; ; \; l \le i \le n \} <
       \delta_n$

\smallskip

(**) $(f_1,\ldots,f_n)=(q(u),X_1-v_1(u),\ldots,X_n-v_n(u))$

\medskip

The network ${\cal N}$ has size $(nd\delta L)^{O(1)}$ and depth
$O(n(log_2(nd) +\ell)log_2 \delta)$.

Let $H=\alpha_1X_1+\cdots+\alpha_nX_n$ be the representation of the nonzero
linear form $H \in k[X_1,\ldots,X_n]$ by its coefficients
$\alpha_1,\ldots,\alpha_n \in k$ and let $p$ be the characteristic (or minimal)
polynomial of the $k$--linear endomorphism of $k[T]/(q)$ induced by the
polynomial $\alpha_1v_1+\cdots+\alpha_nv_n \in k[T]$. We can compute the
coefficients of $p$ from the data $\alpha_1,\ldots,\alpha_n$, $q$ and
$v_1,\ldots,v_n$ by a well parallelizable algorithm in sequential time
$\delta_n^{O(1)}$. Thus we may suppose without loss of generality that ${\cal
N}$ computes also the coefficients of the monic polynomial $p$. The condition
(**) above expresses an isomorphism between the $k$--algebras $k[T]/(q)$ and
$k[X_1,\ldots,X_n]/(f_1,\ldots,f_n)$. From the particular form of this
isomorphism we infer that $p$ is also the characteristic (or minimal) polynomial
of the $k$--linear endomorphism of
$k[X_1,\ldots,X_n]/(f_1,\ldots,f_n)$ induced by the multiplication by $H$. This
implies that $p(H) \in (f_1,\ldots,f_n)$ holds, whence the conclusion of Theorem
\ref{afin}.

\medskip

Theorem \ref{tor} is deduced from Theorem \ref{parametrizacion} $(ii)$ in almost
textually the same manner.
\cqfd

Lemma \ref{sinfor} and Proposition \ref{inductive} can be combined in the same
way as in the proof of Theorem \ref{parametrizacion} $(i)$ in order to obtain
the following elimination result for reduced complete intersection ideals of
positive dimension~:

\begin{proposition}
\label{inter.comp.}

Let 
$f_1,\ldots,f_s$ be polynomials which belong to $k[X_1,\ldots,X_n]$. Suppose
that  $f_1,\ldots,f_s$ form a regular sequence in $k[X_1,\ldots,X_n]$. For any
$1 \le j \le s$ let $\delta_j\!:=deg V(f_1,\ldots,f_j)$ be the geometric degree
of the affine variety defined by the ideal $(f_1,\ldots,f_j)$ which we assume to
be radical. Write $r:=n-s$, $\delta:= max \{\delta_j \; ; \; 1\le j  \le s\}$
and $d:= max \{deg f_j \; ; \; 1 \le j \le s \}$. Suppose that the polynomials
$f_1,\ldots,f_s$ are given by a division free straight--line program $\beta$ in
$k[X_1,\ldots,X_n]$ of length $L$ and depth $\ell$. Then there exists an
arithmetic network with parameters in ${\bar k}$ which has size
$(sd \delta L)^{O(1)}$ and depth
$O(s(log_2 (sd)+\ell)log_2\delta)$ and which produces from the circuit $\beta$
as input, the following items~:

\begin{itemize}

\item a nonsingular matrix of $k^{n \times n}$ which is given by its
      coefficients and which transforms the variables $X_1,\ldots,X_n$ into
      new ones, say $Y_1,\ldots, Y_n$,

\item a nonzero $k$--linear form $u \in k[Y_{r+1},\ldots,Y_n]$,

\item a division free straight--line program $\gamma$ in
      ${\bar k}[Y_1,\ldots,Y_r]$ which represents a nonzero polynomial 
      $\rho \in k[Y_1,\ldots,Y_r]$ and the coefficients with respect to the
      variable $u$ of certain polynomials $q, v_1,\ldots,v_n \in   
k[Y_1,\ldots,Y_r,u]$.
\end{itemize}

These items have the following properties~:

\begin{enumerate}

\item the new variables $Y_1,\ldots, Y_n$ are in Noether position with 
      respect to the ideal $(f_1,$ $\ldots,f_s)$, the variables $Y_1,\ldots, Y_r$
      being free,
\item the polynomial $q$ is monic and separable in $u$ and $\rho$ is its
      discriminant,
\item in the localized ring $k[X_1,\ldots,X_n]_{\rho}$
      the following ideals are identical
$$(f_1,\ldots,f_s)_{\rho}=(q,\rho X_1-v_1,\ldots,\rho
      X_n-v_n)_{\rho},$$
\item the polynomials $q$, $\rho$ and $v_1,\ldots,v_n$ satisfy the degree
      conditions  $deg q= deg_u q= \delta_s$, $deg \rho \le \delta_s^3$, 
$max\{deg_u v_j \; ; \; 1 \le j
      \le n\} < \delta_s$ and
      $max \{deg v_j \; ; \; 1 \le j \le n \}\le 2\delta_s^3$,
\item the length of the  straight--line program $\gamma$ is
      $O((s^5+L)\delta_s^{11})=(sL\delta)^{O(1)}$ and its depth is
      $O((log_2s+\ell)log_2\delta_s)=O((log_2s+\ell)log_2\delta)$.

\end{enumerate} 
\end{proposition}

We omit the proof of this proposition which is essentially the same as the proof
of Theorem \ref{parametrizacion} $(i)$.

\medskip

Of course Proposition \ref{inter.comp.} admits also a toric version with almost
the same proof as Theorem \ref{parametrizacion} $(ii)$. We are going now to
formulate a slightly more general statement in which the r\^ole of the form
$\prod\limits_{i=1}^n X_i$ is played by an arbitrary nonzero polynomial $g
\in k[X_1,\ldots,X_n]$. We think this statement is interesting by its own
because it says that our algorithms are able to ``avoid" undesired points which
are contained in a previously given hypersurface $V(g)$.

Let be given nonzero polynomials $f_1,\ldots,f_s,g \in k[X_1,\ldots,X_n]$ with
$s \le n$. We say that $f_1,\ldots,f_s$ form a {\it secant family (``suite
s\'ecante") which avoids the hypersurface $V(g)$} if for any $ 1 \le j
\le s$ the localization $V(f_1,\ldots,f_j)_g$ is not empty and if for any
irreducible component $C$ of $V(f_1,\ldots,f_j)$ the following implication
holds~:

{\it  if $C$ is {\it not} contained in the hypersurface $V(g)$ then $dimC=n-j$}.
\medskip

Let $1 \le j \le s$. We say that an irreducible component $C$ of
$V(f_1,\ldots,f_j)$ {\it avoids the hypersurface $V(g)$} if $C$ is not contained
in $V(g)$. The union of all irreducible components of
$V(f_1,\ldots,f_j)$ which avoid the hypersurface $V(g)$ is called the {\it part
of the affine variety $V(f_1,\ldots,f_j)$ which avoids the hypersurface $V(g)$}
and we denote this part by $V'(f_1,\ldots,f_j)$. We write
$deg'V(f_1,\ldots,f_j)$ for its geometric degree. Thus $V'(f_1,\ldots,f_j)$ is
the Zariski closure of the locally closed subset $V(f_1,\ldots,f_j)_g$ of ${\bf
A}^n$ and we have
$deg'V(f_1,\ldots,f_j)=degV'(f_1,\ldots,f_j)$. Moreover if $f_1,\ldots,f_s$ form
a secant family which avoids the hypersurface $V(g)$ then
$deg'V(f_1,\ldots,f_j)>0$ holds for any $1 \le j \le s$.

With these notions and notations we are now able to formulate our statement
which generalizes Theorem \ref{parametrizacion} $(ii)$ to positive dimension~:

\begin{proposition}
\label{hipersurf} Let $f_1,\ldots,f_s$ and $g$ be polynomials of
$k[X_1,\ldots,X_n]$ such that
$f_1,\ldots,f_s$ form a secant family which avoids the hypersurface $V(g)$. For
any $1 \le j \le s$ let 
$\delta'_j:=deg'V(f_1,\ldots,f_j)$ be the degree of the part of the affine
variety
$V(f_1,\ldots,f_j)$ which avoids the hypersurface $V(g)$ and suppose that the
localized ideal $(f_1,\ldots,f_j)_g$ is radical. Write
$r:= n-s$, $\delta':=max \{ \delta'_j \; ; \; 1 \le j \le s \}$ and
$d:=max \{deg f_j  \; ; \; 1 \le j \le s \}$. Assume that the polynomials
$f_1,\ldots,f_s$ and $g$ are given by a division free straight--line program
$\beta$ in $k[X_1,\ldots,X_n]$ of length $L$ and depth $\ell$. Then there exists
an arithmetic network with parameters in ${\bar k}$ which has size $(sd
\delta' L)^{O(1)}$ and depth $O(s(log sd+ \ell)log_2
\delta')$ and which produces from the circuit $\beta$ as input, the following
items~:

\begin{itemize}

\item a nonsingular matrix of $k^{n \times n}$ which is given by its
      coefficients and which transforms the variables $X_1,\ldots,X_n$ into new
      ones, say $Y_1,\ldots, Y_n$,

\item a nonzero $k$--linear form $u \in k[Y_{r+1},\ldots,Y_n]$,

\item a division free straight--line program $\gamma$ in 
      ${\bar k}[Y_1,\ldots,Y_r,u]$ which represents a nonzero polynomial $\rho
      \in k[Y_1,\ldots,Y_r]$ and the coefficients with respect to $u$ of
      certain polynomials $q, v_1,\ldots,v_n \in k[Y_1,\ldots,Y_r,u]$.
\end{itemize}

These items have the following properties~:

\begin{enumerate}

\item the new variables $Y_1,\ldots, Y_n$ are in Noether position with respect
      to $V'(f_1,\ldots,f_s)$, the variables $Y_1,\ldots,Y_r$ being free,

\item the polynomial $q$ is monic and separable in $u$ and $\rho$ is its
      discriminant,

\item in the localized ring $k[X_1,\ldots,X_n]_{\rho \cdot g}$ the following
ideals are identical
     $$(f_1,\ldots,f_s)_{\rho \cdot g}=(q,\rho X_1-v_1,\ldots,\rho
      X_n-v_n)_{\rho \cdot g},$$

\item the polynomials $q$, $\rho$ and $v_1,\ldots,v_n$ satisfy the degree
      conditions $deg q=deg_u q= \delta '_s$, $deg\rho \le 2\delta_s ^{'3}$,  
$max \{deg_u v_j \; ; \; 1 \le j
      \le n \} < \delta '_s$ and
      $max \{deg v_j \; ; \; 1 \le j \le n \} \le 2\delta_s ^{'3}$ 

\item the length of the straight--line program $\gamma$ is
      $O((s^5+L)\delta_s^{'11})=(sL \delta')^{O(1)}$ and its depth is 
      $O((log_2 s+\ell)log_2 \delta'_s)=O((log_2 s+\ell)log_2 \delta')$.

\end{enumerate} 
\end{proposition}

We omit the proof of this proposition which essentially the same as the proof of
Theorem \ref{parametrizacion} $(ii)$ with the form $\prod\limits_{i=1}^n X_i$
replaced by the polynomial $g$.

\section{Division Theorems}
\label{division}

In this section we explain how Proposition \ref{inter.comp.} can be applied in
order to obtain refined complexity and degree bounds in the division theorems of
\cite{figismi}, \cite{krick-pardo1} and \cite{saso} (see also \cite{gihesa} and
\cite{krick-pardo:CRAS}). Let us first state our results.

\begin{theorem}
\label{division1} Let $f_1,\ldots,f_s$ and $g$ be polynomials of
$k[X_1,\ldots,X_n]$ such that
$f_1,\ldots,f_s$ form a regular sequence and $g$ belongs to the ideal
$(f_1,\ldots,f_s)$. For $1 \le j\le n$ denote by $\delta_j\!:=deg
V(f_1,\ldots,f_j)$ the geometric degree of the affine variety defined by the
ideal $(f_1,\ldots,f_j)$ which we suppose to be radical. Write 
$\delta \!:= max \{\delta_j \; ; \; 1 \le j \le s-1 \}$ and 
$d\!:= max \{degf_j \; ; \; 1 \le j \le n \}$. Assume that the polynomials
$f_1,\ldots,f_s$ and $g$ are given by a division free straight--line program
$\beta$ in $k[X_1,\ldots,X_n]$ of length $L$ and depth $\ell$. Then there exists
an arithmetic network with parameters in ${\bar k}$ which has size $(sd \delta
L)^{O(1)}$ and depth $O(s(log_2 sd +\ell)log_2\delta)$ and which produces from
the circuit $\beta$ as input a division free straight--line program $\gamma$ in
${\bar k}[X_1,\ldots,X_n]$ such that $\gamma$ represents certain polynomials
$p_1,\ldots,p_s \in k[X_1,\ldots,X_n]$ with the following properties~:

\begin{enumerate}

\item $g=p_1 f_1+ \cdots +p_s f_s$

\item $max \{deg p_j \; ; \; 1 \le j \le s \} \le (2s^2d+ max\{d,degg\}) \delta$
  
\end{enumerate}

Moreover the circuit $\gamma$ has size $(deg g)^2 (sd \delta L)^{O(1)}$ and
depth  $O(s(log_2 sd +\ell)log_2\delta)$. Let us remark that allowing {\it
divisions} in the circuit $\gamma$ we may diminuish its size to $(sd \delta
L)^{O(1)}$.
\end{theorem}

\medskip 

Let $f_1,\ldots,f_s$ be nonconstant polynomials representing an equation system
of  $k[X_1,$ $\ldots,X_n]$. Define polynomials $g_1,\ldots,g_n \in
k[X_1,\ldots,X_n]$ as follows~: if the characteristic of $k$ is zero, choose
$g_1,\ldots,g_n$ as generic linear combinations of the polynomials
$f_1,\ldots,f_s$ and if the characteristic is positive, choose $g_1,\ldots,g_n$
as generic linear combinations of the set of polynomials $\{ X_i f_j \; ; \; 1
\le i \le n\; , \; 1 \le j \le s \}$. For $1\le i\le n$ denote by
$\delta_i\!:=deg V(g_1,\ldots,g_n)$ the geometric degree of the affine variety
$V(g_1,\ldots,g_i)$. Then we call $\delta\!:=max\{\delta_i \; ; \; 1 \le i \le n
\}$ the {\it generic geometric degree} of the equation system given by
$f_1,\ldots,f_s$. Observe that the generic geometric degree of the system
$f_1,\ldots,f_s$ is positive even if $f_1,\ldots,f_s$ generate the trivial ideal.

\medskip

With these notion we can state our next result, namely

\begin{theorem}[effective Nullstellensatz]
\label{division2} Let $f_1,\ldots,f_s$ be polynomials of the ring
\\ $k[X_1,\ldots,X_n]$ which generate the  trivial ideal. Let $\delta$ be the
generic geometric degree of the equation system $f_1,\ldots,f_s$ and let $d:= max
\{deg f_j \; ;
\; 1 \le j \le s \}$. Suppose that the polynomials $f_1,\ldots,f_s$ are given by
a division free straight--line program 
$\beta$ in $k[X_1,\ldots,X_n]$ of length $L$ and depth $\ell$. Then there exists
an arithmetic network with parameters in ${\bar k}$ which has size $(nd \delta
L)^{O(1)}$ and depth $O(n(log nd + \ell)log_2\delta)$ and which produces from
the circuit $\beta$ as input a division free straight--line program $\gamma$ in
${\bar k}[X_1,\ldots,X_n]$ such that $\gamma$ represent certain polynomials 
$p_1,\ldots,p_s
\in k[X_1,\ldots,X_n]$ which have the following properties~:

\begin{enumerate}

\item $1=p_1 f_1+ \cdots +p_s f_s$

\item $max \{deg p_j \; ; \; 1 \le j \le s \} \le 3n^2 d \delta$
\end{enumerate}

Moreover the circuit $\gamma$ has size $(nd \delta L)^{O(1)}$ and depth $O(n(log
nd + \ell)log_2\delta)$.
\end{theorem}

If we replace in the statement of Theorem \ref{division1} and \ref{division2}
the parameter $\delta$ by the B\'ezout estimations $d^s$ and $d^n$ respectively
we obtain the worst case complexity bounds of \cite{figismi}, Th\'eor\`eme 5.1
and 5.2, \cite{krick-pardo1} Theorem 1 and 2 and \cite{saso}, Lemma 15 and 16,
Theorem 19 and 24.

Let us now make some comments first on the degree and then on the complexity
bounds of Theorem \ref{division1} and \ref{division2} and their proofs. The
degree bounds in Theorem \ref{division1} $(ii)$ and \ref{division2} $(ii)$ are
noteworth by their own. They are due to the paper \cite{saso} from which they
follow immediately by performing the following modification in the proof of
\cite{saso}, Lemma 15 and 16~:

\noindent every time when the B\'ezout inequality (\cite{joos:tesis}, Theorem 1)
is applied in order to produce a degree bound of $d^j$, this bound is replaced
by the number $\delta_j$, where $1 \le j\le s$ or $1 \le j \le n$, following the
context (see \cite{saso2}).

We are now going to explain how the complexity bounds of Theorem \ref{division1}
$(i)$ and \ref{division2} $(i)$ follow from Proposition \ref{inter.comp.}. In
order to make explanations easier we refer only to proofs in paper
\cite{figismi} and to the sequential complexity model (compare with
\cite{krick-pardo1} for considerations of nonscalar parallel complexity). We
obtain a proof of Theorem \ref{division1} $(i)$ and \ref{division2} $(i)$ almost
directly from the proof of \cite{figismi} Th\'eor\`eme 5.1. and 5.2. if we
perform systematically the following modifications in the paper \cite{figismi}
~: Proposition 2.3.1. and Proposition principale 2.4.1. are replaced everywhere
they are used in \cite{figismi} by our Proposition \ref{inter.comp.}. Moreover
in all degree estimations, when the application of the B\'ezout inequality
produces a $d^j$ bound (for $1\le j \le s$ or $1 \le j \le n$), this bound is
replaced by $\delta_j$. With these modifications the proof of our Theorem
\ref{division1} $(i)$ and \ref{division2} $(i)$ is almost textually the same as
that of \cite{figismi}, Th\'eor\`eme 5.1 and 5.2. There is only one point where
caution is necessary~: the application of duality techniques in \cite{figismi}
is based on the simple minded decomposition formula for duality $(3)$ in Section
$3.2$ of the same paper. The size of this decomposition may happen to be too big
for our purpose here. In order to remedy this problem we have to use a different
view of duality which is close to \cite{saso} and \cite{saso2}. This will be
contained in a forthcoming paper where complete and self contained proofs of
Theorem \ref{division1} and \ref{division2} will be given.

For getting simultaneously with the complexity bounds also the degree
estimations in Theorem \ref{division1} and \ref{division2} one has to adapt the
algorithm underlying the proof of Theorem \ref{division1} $(i)$ and
\ref{division2} $(i)$ to the particular constructions in \cite{saso} and 
\cite{saso2}.

\section{Conclusions}
\label{conclusiones}

This final section is devoted to the announcement of further complexity results
which can be obtained by our method. We hope that these results, which are
intermediary products of work still in progress, shed light on the possibilities
and limitations of our approach.

First let us turn back to our discussion at the end of Section \ref{compresion}.

Remark \ref{fact_gates} and Proposition \ref{inductive} imply that our main
results, namely Theorem \ref{afin}, \ref{tor}, \ref{parametrizacion},
\ref{division1}, \ref{division2} and Proposition \ref{inter.comp.},
\ref{hipersurf} remain true if we replace in their conclusions ``arithmetic
network with parameters in ${\bar k}$" by ``arithmetic network over $k$ with
factorization gates" and ``straight--line program in ${\bar k}[X_1,\ldots,X_n]$"
by ``straight--line program in $k[X_1,\ldots,X_n]$". Parallel to this
observation we may ask what happens to our algorithms when we try to transfer
them to the bit complexity model (here we suppose that $k$ is $\Q$ or a finitely
generated transcendental extension of $\Q$). We announce here that our
procedures allow perfect a modularization modulo suitably (randomly) chosen
primes of not too big height. This implies that a suitably randomized version of
our main complexity results remain valid in the bit model if we take into account
the (bit) size of the rational parameters which appear in the input circuit
$\beta$. Finally let us state just one application of our method to an
elimination problem in semialgebraic geometry.

Let $k$ be an ordered hilbertian field with factorization at moderate costs and
let $K$ be a real closure of $k$. Let ${\bar k}\!:={\bar K}\!:=K(i)$ with
$i^2=-1$ be the corresponding algebraic closure of $k$ and $K$. We consider the
affine space ${\bf A}^n\!:={\bf A}^n({\bar k})$ equipped with the Zariski
topology whose closed sets are the $k$--definable algebraic subsets of ${\bar
k}^n$. Let $W\subseteq {\bf A}^n$ be  a closed subset and let
$W=C_1\cup\cdots \cup C_s$ be its decomposition in irreducible components with
respect to this topology. Thus $W$ and $C_1,\ldots,C_s$ are $k-$ definable
algebraic subsets of  ${\bar k}^n$. Let $1\le j \le s$ and consider the
irreducible component  $C_j$ of $W$. We call $C_j$ a {\it real component of $W$}
if the real variety $C_j\cap K^n$ contains a smooth point of $C_j$. Let 
$$I =\{ j \in \N : 1 \le j \le s, {\hbox {\rm where}} \hskip 3pt C_j \hskip
3pt{\hbox {\rm  is a real component of $W$}} \}.$$ We call the affine variety
$W''\!:=\bigcup \limits_{j \in I} C_j$ {\it the real part of W} and we define
the {\it real degree of W} as
$deg''W\!:=degW''=\sum\limits_{j \in I}deg C_j$. Observe that $deg''W=0$ holds
if and only if the real part $W''$ of $W$ is empty.

Let $f$ be a nonconstant and squarefree polynomial of $k[X_1,\ldots,X_n]$ which
defines a hypersurface $W\!:=V(f)$ of ${\bf A}^n$. Let $V:=W\cap K^n$ the real
variety given by the polynomial $f$. {\it Suppose that $V$ is nonempty, bounded
and smooth with regular equation $f$. Assume that the variables $X_1,\ldots,X_n$
are in  generic position with respect to $V$ and consider for any $0\leq i<n$
the polar variety of $W$ corresponding to the linear space defined by the
equations $X_1,\ldots,X_i$. Denote by $W_i$ this polar variety and consider the
real variety $V_i:=W_i\cap K^n$. With these notions and notations it is not too
difficult to deduce from the Weak Transversality Theorem (\`a la Sard--Thom)
that for any $0\leq i<n$ the following facts are true~:

\begin{itemize}

\item $W_0=W$,

\item $W_i$ is a nonempty equidimensional affine subvariety of ${\bf A}^n$
which   is smooth in all of its points which are smooth with respect to $W$, 

\item the real part $W_i''$ of the polar variety $W_i$ coincides with the 
      Zariski closure of $V_i$ in ${\bf A}^n$,
\item $V_i$ is defined by the equations 
      $f,{{\partial f} \over {\partial
      X_1}}, \ldots,{{\partial f}
      \over {\partial X_i}}$,

\item for any $i<j \le n$ the ideal  
      $$\left(f, {{\partial f} \over {\partial X_1}}, \ldots,{{\partial f}
      \over {\partial X_i}}\right)_{{{\partial f} \over {\partial X_j}}}$$ is
       radical.
\end{itemize}}

Let us write $\delta''_i\!:= deg''W_i$ for $0\le i < n$ and
$\delta''\!:=max\{\delta''_i \; ;\; 0\le i < n\}$ and $d\!:=deg f$.

With these notations and assumptions we have the following complexity result~:

\begin{theorem}(in collaboration with B. Bank and R. Mandel)
\label{real1} Suppose that the equation $f$ is given by a division free
straight--line program $\beta$ in $k[X_1,\ldots,X_n]$ of length $L$ and depth
$\ell$. Then there exists an arithmetic network in $k$ with factorization gates
which has size $(nd\delta'' L)^{O(1)}$ and depth $O(n(log_2(nd)
+\ell)log_2\delta'')$ which from the circuit $\beta$ as input produces the
coefficients of a nonzero linear form $u\in k[X_1,\ldots,X_n]$ and of certain
polynomials $q,v_1,\ldots,v_n \in k[u]$. These polynomials have the following
properties~:

\begin{enumerate}
\item $q$ is monic and squarefree and has degree $\delta''_n$. Moreover
$v_1,\ldots,v_n$ satisfy the degree bound $max\{deg v_j \; ;\; 1\le j < n\} <
\delta''_n$ 
\item For each semialgebraically connected component $C$ of $V$ there exists a
      point $\xi\in C$ and a field element 
      $\tau \in K$ such that $q(\tau)=0$ and $\xi=(v_1(\tau),\ldots,v_n(\tau))$
      holds.
\end{enumerate}
\end{theorem}
 We observe here that the network computes only points $\xi$ which are  critical
points of the projection map of $V$ into $K$ induced by the coordinate
$X_n$.

\bigskip

{\bf Acknowledgments}~: We thank the unknown referee for his inspiring comments
and his suggestions to improve the presentation of this paper.

\end{document}